\documentclass[11pt,draft,reqno]{amsart}
\input epsf.tex

\usepackage[all]{xy}
\usepackage{amsthm,array,amssymb,amscd,amsfonts,latexsym}
\xyoption{poly}

\renewcommand{\baselinestretch}{1.1}

\DeclareFontFamily{OT1}{wncyr}{\hyphenchar\font45
}
\DeclareFontShape{OT1}{wncyr}{m}{n}{%
   <5> <6> <7> <8> <9> gen * wncyr
   <10> <10.95> <12> <14.4> <17.28> <20.74>  <24.88>wncyr10}{}
\DeclareFontShape{OT1}{wncyr}{m}{it}{%
   <5> <6> <7> <8> <9> gen * wncyi
   <10> <10.95> <12> <14.4> <17.28> <20.74> <24.88> wncyi10}{}
\DeclareFontShape{OT1}{wncyr}{m}{sc}{%
   <5> <6> <7> <8> <9> <10> <10.95> <12> <14.4>
   <17.28> <20.74> <24.88>wncysc10}{}
\DeclareFontShape{OT1}{wncyr}{b}{n}{%
   <5> <6> <7> <8> <9> gen * wncyb
   <10> <10.95> <12> <14.4> <17.28> <20.74> <24.88>wncyb10}{}
\input cyracc.def
\def\rus{\usefont{OT1}{wncyr}{m}{n}\cyracc\fontsize{8}{11pt}\selectfont}

\DeclareMathSizes{9}{9}{7}{5} 

\begin{document}

\newtheorem
{theorem}{Theorem}

\newtheorem{proposition}
{Proposition}
\newtheorem{lemma}
{Lemma}
\newtheorem{corollary}
{Corollary}
\newtheorem{question}
{Proposition}

\theoremstyle{definition}
\newtheorem{definition}
{Definition}
\newtheorem{remark}
{Remark}
\newtheorem{example}
{Example}
\newtheorem{reduction}
{Reduction}
\newtheorem{problem}
{Problem}

\theoremstyle{definition}
\newtheorem{enavant}[proposition]{\hskip -.0mm}


\newcommand{\XG}{X_{\lambda_1,\ldots,\lambda_d}{}_{{}^{\overset{}
{\centerdot}}}
  \hskip .5mm G\hskip
-3.9mm^{\overset{\centerdot}{}} \hskip
-1.37mm{}^{\underset{\centerdot}{}}\hskip
3.5mm{} }
\newcommand{\XC}{X_{\lambda_1,\ldots,\lambda_d}/\!\!/G}


\title[Tensor products and  open orbits]{Tensor
product
decompositions and \\
open orbits in multiple flag varieties}


\author[Vladimir  L. Popov]{Vladimir  L. Popov${}^*$}
\address{Steklov Mathematical Institute,
Russian Academy of Sciences, Gubkina 8,
Moscow\\ 119991, Russia}
\email{popovvl@orc.ru}

\thanks{
 ${}^*$\,Supported by ETH, Z\"urich,
 Switzerland, Russian grants {\rus RFFI
05--01--00455}, {\rus N{SH}--9969.2006.1},
and program {\it Contemporary Problems of
Theoretical Mathematics} of the
Mathe\-matics Branch of the Russian Academy
of Sciences.}

\date{August 25, 2006}



\begin{abstract} For a connected
semisimple algebraic group $G$, we consider
some special infinite series of tensor
products
 of simple $G$-modules whose
 $G$-fixed point spaces are at most
one-dimensional. We prove that their
existence is closely related to the
existence of open $G$-orbits in multiple
flag varieties and address the problem of
classifying
 such series.
\end{abstract}

\maketitle

\vskip 10mm

Let $G$ be a connected simply connected
semisimple algebraic group. In this paper
we establish a close interrelation between
some special series of tensor products of
simple $G$-modules whose $G$-fixed point
spaces are at most one-dimensional and
multiple flag varieties of $G$ that contain
open $G$-orbit. Motivated by this intimate
connection with geometry, we then address
the
  problem of classifying such series.
Starting with the basic definition
 and examples in Sections
 \ref{basic} and \ref{examples},
 we  introduce necessary notation in Section \ref{no}
 and then formulate
 our main
results in Section \ref{main}. Other
results and proofs
 are contained in the remaining part of paper.

Below all algebraic varieties are taken
over an algebraically closed field $k$ of
characteristic zero.


\begin{enavant}{\bf Basic definition}
\label{basic}

 Fix a choice of
Borel subgroup $B$ of $G$ and maximal torus
$T\subset B$. Let ${\rm P}_{++}$ be the
additive monoid of dominant characters of
$T$ with respect to $B$. Put $${\rm
P}_{\gg}:={\rm P}_{++}\setminus \{0\}.$$
For $\lambda\in {\rm P}_{++}$, denote by
$E_\lambda$ a simple $G$-module of highest
weight $\lambda$ and by $\lambda^*$ the
highest weight of dual $G$-module
$E_\lambda^{\,*}$. Let $P_\lambda$ be the
$G$-stabilizer of  unique $B$-stable line
in $E_\lambda$. If $\mu, \lambda_1,\ldots,
\lambda_d\in {\rm P}_{++}$, denote by
$c_{\lambda_1,\ldots,\lambda_d}^{\mu}$ the
multiplicity of $E_{\mu}$ inside
$E_{\lambda_1}\otimes\ldots\otimes
E_{\lambda_d}$, i.e., the
Littlewood--Richardson coefficient $\dim
\bigl({\rm Hom}_G(E_\mu,
E_{\lambda_1}\otimes\ldots\otimes
E_{\lambda_d})\bigr)$. Denote respectively
by $\varpi_1,\ldots,\varpi_r$ and
$\alpha_1,\ldots,\alpha_r$
 the systems of fundamental weights of ${\rm
P}_{++}$ and simple roots of $G$ with
respect to $T$ and $B$
 enumerated as in
 \cite{bourbaki}.
Let respectively $\mathbb Z^{}_{\geqslant
0}$ and~$\mathbb Z^{}_{> 0}$ be the sets of
all nonnegative and all positive integers.
We write ${\rm P}_{\gg}^d$ in place of
$({\rm P}_{\gg})^d$, etc.

\begin{definition}\label{simple}
We call a $d$-tuple
$(\lambda_1,\ldots,\lambda_d)\in {\rm
P}_{\gg}^{d}$ {\it primitive} if
\begin{equation}\label{l}\textstyle
c_{n_1\lambda_1,\ldots,
n_d\lambda_d}^0\leqslant 1 \hskip
3mm\mbox{for all}\hskip 2mm
(n_1,\ldots,n_d)\in \mathbb Z_{\geqslant
0}^d.
\end{equation}
\end{definition}

Schur's lemma and the isomorphism
$(E_\mu\otimes E_\nu)^G\simeq {\rm
Hom}^{}_G(E_{\mu^*}, E_\nu)$ imply that
\begin{equation}\label{multipl}
c_{\mu, \nu}^0=\begin{cases} 1&\mbox{if
$\mu^*=\nu$},\\
0&\mbox{otherwise;}
\end{cases}
\end{equation}
whence condition \eqref{l} is equivalent to
the following:
\begin{gather}\label{simple2}
\begin{gathered}
c^{n_j\lambda_j^*}_{n_1\lambda_1,\ldots,
\widehat{n_j\lambda_j},\ldots,
n_d\lambda_d}\hskip -1mm\leqslant 1 \mbox{
for all
 $(n_1,\ldots,n_d)\in \mathbb
Z_{\geqslant 0}^d$ and some (equivalently,
every) $\lambda_j$.}
\end{gathered}
\end{gather}

The set of primitive elements of ${\rm
P}_{\gg}^d$ is clearly stable with respect
to permutation of coordinates and
automorphisms of ${\rm P}_{\gg}^d$ induced
by automorphisms of the Dynkin diagram of
~$G$. If $(\lambda_1,\ldots, \lambda_d)\in
{\rm P}_{\gg}^d$ is primitive, then
$(\lambda_{i_1},\ldots, \lambda_{i_s})\in
{\rm P}_{\gg}^s$ is primitive for every
subset $\{i_1,\ldots, i_s\}$ of
$\{1,\ldots, d\}$.

\begin{remark}\label{r1} The notion of primitive
$d$-tuple admits a natural generalization:
we call a $d$-tuple
$(\lambda_1,\ldots,\lambda_d)\in {\rm
P}_{\gg}^d$ {\it primitive at} $\mu\in {\rm
P}_{++}$ if
\begin{equation}\label{primat}
c_{n_1\lambda_1,\ldots,
n_d\lambda_d}^{\mu}\leqslant 1\hskip 3mm
\mbox{for all \hskip 2mm
$(n_1,\ldots,n_d)\in \mathbb Z_{\geqslant
0}^d$}.
\end{equation}
Then ``primitive'' means ``primitive at
$0$''.
\end{remark}
\end{enavant}


\begin{enavant}
{\bf Examples} \label{examples}

Clearly, for $d=1$, every element of ${\rm
P}_{\gg}^d$ is primitive. By
\eqref{multipl} the same is true for $d=2$.
For $d\geqslant 3$, the existence of
primitive elements in ${\rm P}_{\gg}^d$ is
less evident.

\begin{example}\label{e1} Let $G={\bf
SL}_2$. Then ${\rm P}_{++}=\mathbb
Z^{}_{\geqslant 0} \varpi_1$ and
$G/P_{\lambda}$ for $\lambda\neq 0$ is
isomorphic to the projective line ${\bf
P}^1$. From Definition~\ref{simple} and the
Clebsch--Gordan formula
\begin{equation*}\textstyle
E_{s\varpi_1}\otimes E_{t\varpi_1}\simeq
\bigoplus_{0\leqslant i\leqslant t}
E_{(s+t-2i)\varpi_1},\quad s\geqslant t,
\end{equation*}
it is not difficult to deduce that an
element of ${\rm P}_{\gg}^d$ is primitive
if and only if $d\leqslant 3$.
Theorems~\ref{teo1} and \ref{teo2} below
imply that this is equivalent to the
classical fact that for the diagonal action
of ${\bf SL}_2$
 on $({\bf P}^1)^d$, an open orbit
 exists if and only if
$d\leqslant 3$.
\end{example}

\begin{example} \label{e2} If, for
$(\lambda_1,\ldots,\lambda_d)\in {\rm
P}_{\gg}^d$, condition \eqref{primat} holds
for every $\mu\in {\rm P}_{++}$, then by
\eqref{simple2} the $(d+1)$-tuple
 $(\nu,
\lambda_1,\ldots,\lambda_d)$ is primitive
for every $\nu\in {\rm P}_{\gg}$. Such
$d$-tuples $(\lambda_1,\ldots,\lambda_d)$
exist. For instance, if $G$ is of type {\sf
A}, {\sf B}, {\sf C}, {\sf D}, or ${\sf
E}_6$, then the explicit decomposition
formulas for $E_{m_1\varpi_1}\otimes
E_{m_2\varpi_1}$ (see
\cite[1.3]{littelmann} or, for the types
{\sf A}, {\sf B}, {\sf C}, {\sf D},
\cite[pp.\,300--302]{OV}) imply that
$(\varpi_1, \varpi_1)$ shares this
property. For $G={\bf SL}_n$, the
classification of all $d$-tuples
$(\lambda_1,\ldots,\lambda_d)$ sharing this
property can be deduced from \cite{St},
where the classification of all
multiplicity free tensor products of simple
${\bf SL}_n$-modules is obtained.
\end{example}

\begin{example}\label{e3}
Let $G$ be the group of type ${\sf E}_{6}$.
By \cite[1.3]{littelmann}, for every $s, t
\in\mathbb Z^{}_{\geqslant 0}$, the
following decomposition holds:
\begin{equation}\label{tensor5}
E_{s\varpi_1}\otimes
E_{t\varpi_1}\simeq\hskip
-7mm\bigoplus_{\fontsize{8pt}{2mm}
\selectfont\biggl\{\begin{gathered}
a_1,\ldots, a_4\in\mathbb Z^{}_{\geqslant 0}\\[-2.5pt]
a_1+a_3+a_4=s\\[-2.5pt]
a_2+a_3+a_4=t\end{gathered} }\hskip -7mm
E_{(a_1+a_2)\varpi_1+a_3\varpi_3+
a_4\varpi_6}.
\end{equation}
Since
$\bigl((a_1+a_2)\varpi_1+a_3\varpi_3+a_4\varpi_6\bigr)^*=
a_4\varpi_1+a_3\varpi_5+(a_1+a_2)\varpi_6$,
it follows from \eqref{tensor5} and
\eqref{multipl} that $\dim
\bigl(\bigotimes_{1\leqslant i\leqslant
4}E_{n_i\varpi_i} \bigr)\hskip 0mm^G$ is
equal to the number of solutions in
$\mathbb Z^{}_{\geqslant 0}$
 of the following system
of eight linear equations in eight
variables $a_1,\ldots, a_4, b_1,\ldots,
b_4$:
\begin{gather}
\left\{ \hskip -57mm
\begin{align*}
a_4&=b_1+b_2,\\[-4.3pt]
a_3&=0,\\[-4.3pt]
b_3&=0,\\[-4.3pt]
a_1&+a_2=b_4,\\[-4.3pt]
a_1&+a_3+a_4=n_1,\\[-4.3pt]
a_2&+a_3+a_4=n_2,\\[-4.3pt]
b_1&+b_3+b_4=n_3,\\[-4.3pt]
b_2&+b_3+b_4=n_4.
\end{align*}
\right.
\end{gather}
Since this system is nondegenerate, it has
at most one such solution. Thus for $G$ of
type ${\sf E}_{6}$, the $4$-tuple
$(\varpi_1, \varpi_1, \varpi_1, \varpi_1)$
is primitive. By Theo\-rems~\ref{teo1} and
\ref{teo2} below (see also
\cite[
Theorem~6]{popov2}) this is equivalent to the
existence of an open $G$-orbit in
$(G/P_{\varpi_1})^4$. Observe that this example in
not in the range of Example~\ref{e2}: for instance,
$c_{4\varpi_1, 4\varpi_1,
3\varpi_1}^{\varpi_1+3\varpi_3+\varpi_5}=2$ (this may
be verified, e.g., utilizing {\sf LiE}).
\end{example}

\begin{example}\label{st-unst} The following
definition singles out a natural subset in
the set of all pri\-mi\-tive $d$-tuples.
Theorem \ref{geom} below shows that this
subset admits a geometric characterization
as well.

\begin{definition}\label{def2}
We call a $d$-tuple
$(\lambda_1,\ldots,\lambda_d)\in {\rm
P}_{\gg}^{d}$ {\it invariant-free} if
\begin{equation*}\label{if}\textstyle
c_{n_1\lambda_1,\ldots, n_d\lambda_d}^0=0
\hskip 3mm\mbox{for all}\hskip 2mm
(n_1,\ldots,n_d)\in \mathbb Z_{\geqslant
0}^d.
\end{equation*}
\end{definition}

Clearly, every $1$-tuple is invariant-free.
For $d=2$, it follows from \eqref{multipl}
that
\begin{equation*}
\mbox{$(\lambda_1, \lambda_2)$ is
invariant-free $\hskip 2mm
\Longleftrightarrow\hskip 2mm \mathbb
Q\,\lambda_1\neq\mathbb Q\,\lambda_2^*$}.
\end{equation*}
\end{example}

\end{enavant}


\begin{enavant}{\bf Notation and
conventions}\label{no}

Below we utilize the following notation,
conventions, and definitions.

$\bullet$ \hskip .5mm $k[Y]$ and $k(Y)$ are
respectively the algebra of regular
functions
 and field of rational functions of an
irreducible algebraic variety $Y$.

\smallskip

$\bullet$ \hskip .5mm  ${\rm Cl(Y)}$ is the
Weil divisor class group of an irreducible
normal variety $Y$. For a nonconstant
function $f\in k(Y)$, the Weil divisor,
divisor of zeros, and divisor of poles of
$f$ are respectively $(f)$, $(f)_0$, and
$(f)_\infty$.

\smallskip

$\bullet$ \hskip .5mm   If $H$ is an
algebraic group, ${\rm Lie}(H)$ and
$\mathcal X(H)$ are respectively the Lie
algebra and the character group ${\rm
Hom}_{\rm alg}(H, {\bf G}_m)$ of $H$. We
utilize additive notation for $\mathcal
X(H)$
 and identify in the natural way $\mathcal
X(H)$  with the lattice in  rational vector
space
$$\mathcal X(H)_{\mathbb Q}:=\mathcal
X(H)\otimes \mathbb Q.$$

If $H$ is a connected reductive group,
$B_H$ its Borel subgroup, and $S\subseteq
B_H$ a maximal torus,  we identify the set
of isomorphism classes of simple algebraic
$H$-modules with a submonoid $\mathcal
X(S)_{++}$ of $\mathcal X(S)$ assigning to
every simple $H$-module $V$ the $S$-weight
of unique $B_H$-stable line in $V$.

\smallskip

$\bullet$ \hskip .5mm  Below all algebraic
group actions
 are algebraic. The
action of $H$ on $H/F$ is that by left
multiplication. If $H$ acts on $Y_1,\ldots,
Y_n$, the action of $H$ on
$Y_1\times\ldots\times Y_n$ is the diagonal
one.

If $H$ acts on a variety $Y$, then $H\cdot
y$ and $H_y$ are respectively the $H$-orbit
and $H$-stabilizer of a point $y\in Y$, and
$k[Y]^H$ and $k(Y)^H$  are the subalgebra
and subfield of $H$-invariant elements in
$k[Y]$ and $k(Y)^H$.

\begin{definition}\label{ample} For an
irreducible variety $Y$, we call the action
of $H$ on $Y$ {\it ample} if $k(Y)^H$ is
algebraic over the field of fractions of
$k[Y]^H$.
\end{definition}

Recall the following definition introduced
in \cite{popov0}.
\begin {definition} \label{stableact}
The action of $H$ on $Y$ is called {\it
stable} if $H$-orbits of points lying off a
proper closed subset of $Y$ are closed in
$Y$.
\end{definition}

$\bullet$ \hskip .5mm The natural action of
$H$ on $k[Y]$ is locally finite. If $H$ is
a connected reductive group, by
$k[Y]_{\lambda}$, where $\lambda\in
\mathcal X(S)_{++}$, we denote
 the
$\lambda$-isotypical component  of
$H$-module $k[Y]$ and put
\begin{equation}\label{monoid}
\mathcal S(H, Y):=\{\lambda \in  \mathcal
X(S)_{++} \mid k[Y]_{\lambda}\neq 0\}.
\end{equation}
The set $\mathcal S(H, Y)$ is a submonoid
of $\mathcal X(S)$ (indeed, $\mathcal S(H,
Y)$ is the set of~all weights of the
natural action of $S$ on
$k[Y]^{B_{\!H}^u}$, where $B_{\!H}^u$ is
the unipotent radical of $B_H$; whence the
claim). If $Y$ is an affine variety, then
the monoid $\mathcal S(H, Y)$ is
finitely
generated (this readily follows from the
fact that in this case $k[Y]^{B_{\!H}^u}$
is a finitely generated $k$-algebra, see,
e.g.,\,\cite[3.14]{popov-vinberg2}).

\smallskip

$\bullet$ \hskip .5mm   If $H$ is a
reductive group and $Y$ is an affine
variety, we denote by
\begin{equation}\label{cat}
\pi^{}_{H, Y}\!: Y\longrightarrow Y/\!\!/H
\end{equation}
the corresponding categorical quotient,
i.e., $Y/\!\!/H$ is an affine variety and
$\pi^{}_{H, Y}$ a dominant (in fact,
surjective) morphism such that $\pi_{H,
Y}^*(k[Y/\!\!/H])=k[Y]^H$,
see,\,e.g.,\,\cite[4.4]{popov-vinberg2}.
Utilizing $\pi_{H, Y}^*$, we identify
$k[Y/\!\!/H]$ with $k[Y]^H$ and, if $Y$ is
irreducible, $k(Y/\!\!/H)$ with the field
of fractions of $k[Y]^H$.

It follows from Rosenlicht's theorem
\cite[
Theorem]{rosenlicht}, see also, e.g.,\,\cite[Cor.\,of
Theorem\,2.3]{popov-vinberg2}, that the action of $H$
on $Y$ is ample if and only if $\dim (Y/\!\!/H)=\dim
(Y)-\underset{y\in Y}{\rm max} \dim (H\cdot y)$
(i.e., in the terminology of \cite[Sect.\,4]{luna},
``$Y/\!\!/H$ {\it a la bonne dimension}''). One can
also prove that if $Y$ is normal, then the action of
$H$ on $Y$ is ample if and only if $k(Y)^H$ is the
field of fractions of $k[Y]^H$.

If the action of $H$ on $Y$ is stable, it
is ample (indeed, since every fiber of
$\pi^{}_{H, Y}$ contains a unique closed
$H$-orbit, see, e.g.,\,\cite[Cor.\,of
Theorem\,4.7]{popov-vinberg2}, the action of
 is stable if and only if every general
fiber   is a closed $H$-orbit of maximal
dimension).

\smallskip

$\bullet$ \hskip .5mm  $N_G(T)$ is the
normalizer of $T$ in $G$ and $W:=N_G(T)/T$
 is the Weyl group of $G$. For every
$w\in W$ we fix a choice of its
representative $\overset{.}w$ in $N_G(T)$.
 The space $\mathcal
X(T)_{\mathbb Q}$ is endowed with the
natural $W$-module structure.

 For
$\lambda=\sum_{i=1}^{r}a_i\varpi_i$,
$a_i\in \mathbb Z^{}_{\geqslant 0}$, the
{\it support} of $\lambda$ is
\begin{equation*}\label{support}
{\rm supp}(\lambda):=\{i\in \mathbb Z^{}_{>
0}\mid a_i\neq 0\}.\end{equation*} The
subgroup $P_\lambda$ is then generated by
$T$ and one-dimensio\-nal unipotent root
subgroups of $G$ corresponding to all
positive roots and those negative roots
that are linear combinations of
$-\alpha_i$'s with $i\notin {\rm
supp}(\lambda)$. We have the equivalence
\begin{equation}\label{parab}
P_\mu=P_\nu \hskip 2mm
\Longleftrightarrow\hskip 2mm {\rm
supp}(\mu)={\rm supp}(\nu).
\end{equation}

For every subset $A$ of $\subseteq {\rm
P}_{++}$, we put
\begin{equation*}
A^*:=\{\lambda^*\mid \lambda\in A\}.
\end{equation*}

$\bullet$\hskip 2mm  We fix a choice of
nonzero point $v_\lambda$ of the unique
$B$-stable line in $E_\lambda$ and denote
by ${\mathcal O}_\lambda$ the $G$-orbit of
$v_\lambda$ and by $\overline{\mathcal
O}_\lambda$ its closure
in $E_\lambda$. We put
\begin{equation}\label{X}
\mathcal O_{\lambda_1,\ldots,\lambda_d}:=
\mathcal O_{\lambda_1}\times\ldots \times
\mathcal O_{\lambda_d} \quad\mbox{and}\quad
X_{\lambda_1,\ldots,\lambda_d}:=\overline{\mathcal
O}_{\lambda_1}\times\ldots\times
\overline{\mathcal O}_{\lambda_d}
\end{equation}
and  identify in the natural way
$X_{\lambda_1,\ldots,\lambda_d}$  with the
closed subset of
$E_{\lambda_1}\oplus\ldots\oplus
E_{\lambda_d}$.

\smallskip

$\bullet$\hskip 2mm If $M$ is a subset of a
vector space, then ${\rm conv}(M)$ and
${\rm cone}(M)$ are respectively the covex
hull of $M$ and the convex cone generated
by $M$. If $M$ is a convex set, ${\rm
int}(M)$ is the set its (relative) interior
points.

\smallskip

$\bullet$\hskip 2mm We put
$$\mathbb Q^{}_{> 0}:=\{a\in \mathbb Q\mid
a> 0\}\quad \mbox{and}\quad \mathbb
Q^{}_{\geqslant 0}:=\{a\in \mathbb Q\mid
a\geqslant 0\}.$$

$\bullet$\hskip 2mm $|N|$ is the
cardinality of a finite set $N$.
\end{enavant}


\begin{enavant}{\bf  Main results}
\label{main}

In this section we formulate main results
of this paper.

Theorem \ref{teo1} explicitly formulates
the aforementioned remarkable connection of
primitive tuples with geometry.

\begin{theorem}\label{teo1}
Let  $(\lambda_1,\ldots, \lambda_d)\in {\rm
P}_{\gg}^d$.
\begin{enumerate}
\item[\rm (i)] If $G/P_{\lambda_1}
\times\ldots\times G/P_{\lambda_d}$
contains an open $G$-orbit, then
$(\lambda_1,\ldots, \lambda_d)$ is
primitive. \item[\rm(ii)] If
$(\lambda_1,\ldots, \lambda_d)$ is
primitive  and the action of $G$ on
$X_{\lambda_1,\ldots, \lambda_d}$ is ample,
then $G/P_{\lambda_1} \times\ldots\times
G/P_{\lambda_d}$ contains an open
$G$-orbit.
\end{enumerate}
\end{theorem}

Theorem \ref{teo1} and equivalence
\eqref{parab} imply

\vskip 2mm

\noindent{\bf Corollary.} \label{suppp}
{\it Let $(\lambda_1,\ldots, \lambda_d)$
and $(\mu_1,\ldots, \mu_d)\in {\rm
P}_{\gg}^d$. Assume that}
\begin{equation*}\label{supp}
{\rm supp}(\lambda_i)={\rm supp}(\mu_i) \
\mbox{\it for all $i$.}
\end{equation*}
{\it If $(\lambda_1,\ldots, \lambda_d)$ is
primitive and the action of $G$ on
$X_{\lambda_1,\ldots, \lambda_d}$ is ample,
then $(\mu_1,\ldots, \mu_d)$  is primitive
as well.}

\vskip 2mm

Theorem \ref{teo1} clarifies relation of
classifying primitive $d$-tuples to the
following problems.

\begin{problem}\label{p1}
Classify multiple flag varieties
$G/P_{\lambda_1}\times\ldots \times
G/P_{\lambda_d}$ that contain open
$G$-orbit.
\end{problem}
\begin{problem} \label{p2}
For what $d$-tuples
$(\lambda_1,\ldots,\lambda_d)\in {\rm
P}_{\gg}^d$ is the action of $G$ on
$X_{\lambda_1,\ldots,\lambda_d}$ ample?
\end{problem}

Regarding Problem \ref{p1},  obvious
dimension reason yields the finiteness
statement about length of possible
$d$-tuples: $d$ in Problem \ref{p1} cannot
exceed a constant depending only on $G$. A
more thorough analysis leads to the
following upper bounds.

\begin{theorem}\label{bounds} Let $G$ be
a simple group. If
$G/P_{\lambda_1}\times\ldots \times
G/P_{\lambda_d}$ contains an open
$G$-orbit, then $d\leqslant b^{}_G$ for the
following $b^{}_G$:

\vskip 3mm

\begin{center}

{\sc Table 1}\vskip 3mm
\begin{tabular}{c||c|c|c|c|c|c|c|c|c}
\text{ {\rm type of} $G$} & ${\sf A}_l$,
\mbox{\fontsize{8pt}{5mm}\selectfont
$l\geqslant 1$}& ${\sf B}_l,$
\mbox{\fontsize{8pt}{5mm}\selectfont
$l\geqslant 3$}& ${\sf C}_l,$
\mbox{\fontsize{8pt}{5mm}\selectfont
$l\geqslant 2$}& ${\sf D}_l$,
\mbox{\fontsize{8pt}{5mm}\selectfont
$l\geqslant 4$}& ${\sf E}_6$ & ${\sf E}_7$&
${\sf E}_8$& ${\sf F}_4$& ${\sf G}_2$
\\[2pt]
 \hline
&&&&&&&&&\\[-9pt]
$b^{}_G$ &$l+2$&$l+1$&$l+1$&$l$&$4$&$4$&$4$&$3$&$2$\\[-9pt]
&&&&&&&&
\end{tabular}\\
\vskip 2mm
\end{center}
\end{theorem}

Notice that, by the Bruhat decomposition,
every multiple flag variety
$G/P_{\lambda_1}\times G/P_{\lambda_2}$
contains only finitely many $G$-orbits (one
of which therefore is open).

In \cite{popov2} a complete solution to
Problem \ref{p1} for $d$-tuples of the form
\begin{equation}\label{mfund}
(\lambda_1,\ldots,\lambda_d)=(m_1\varpi_i,\ldots,
m_d\varpi_i), \quad 1\leqslant i\leqslant
r,\quad (m_1,\ldots, m_d)\in \mathbb Z_{>
0}^d. \end{equation}
 is obtained; the answer is the following.

\begin{theorem}[\cite{popov2}] \label{multi}
Let $G$ be a simple group.  Assume that
$d\geqslant 3$ and
$(\lambda_1,\ldots,\lambda_d)$ is given by
equality \eqref{mfund}. Then the multiple
flag variety
$G/P_{\lambda_1}\times\ldots\times
G/P_{\lambda_d}$ contains an open $G$-orbit
if and only if the following conditions
hold:

\vskip 1mm

\begin{center}

{\sc Table 2} \vskip 3mm
\begin{tabular}{c||c|c|c|c|c|c}
\text{ {\rm type of} $G$} & ${\sf A}_l$,
\mbox{\fontsize{8pt}{5mm}\selectfont
$l\geqslant 1$}& ${\sf B}_l$,
\mbox{\fontsize{8pt}{5mm}\selectfont
$l\geqslant 3$}& ${\sf C}_l$,
\mbox{\fontsize{8pt}{5mm}\selectfont
$l\geqslant 2$}& ${\sf D}_l$,
\mbox{\fontsize{8pt}{5mm}\selectfont
$l\geqslant 4$}& ${\sf E}_6$ & ${\sf E}_7$
\\[2pt]
 \hline
 &&&&&&\\[-9pt]
 {\rm conditions}
 &
 $\begin{matrix}
 \mbox{
 \fontsize{12pt}{10mm} \selectfont
 $d<\frac{(l+1)^2}{i(l+1-i)}$
 }
 \end{matrix}$
 &$\begin{matrix}
 d=3,\\[-2pt] i=1, l
 \end{matrix}$&$
 \begin{matrix}
 d=3,\\[-2pt] i=1, l
 \end{matrix}$&
$\begin{matrix}
 d=3,\\[-2pt] i=1, l-1, l
 \end{matrix}$ & $\begin{matrix}
 d\leqslant 4,\\[-2pt] i=1, 6
 \end{matrix}$ &$\begin{matrix}
 d=3,\\[-2pt] i=7
 \end{matrix}$
\\
\end{tabular}\\
 \vskip 2.5mm
\end{center}
\end{theorem}

The next two theorems concern Problem
\ref{p2}.

\begin{theorem}
\label{teo2} Let
$(\lambda_1,\ldots,\lambda_d)\in {\rm
P}_{\gg}^d$. If every $\lambda_s$ is a
multiple of a fundamental weight,
\begin{equation}\label{mff}
\lambda_s\in \mathbb Z_{>0}\varpi_{i_s}, \
\ s=1,\ldots, d,
\end{equation}
then the action of $G$ on
$X_{\lambda_1,\ldots,\lambda_d}$ is ample.
\end{theorem}

\begin{theorem}\label{fntn} Let $G$
be a simple group and let
$(\lambda_1,\ldots,\lambda_d)\in {\rm
P}_{\gg}^d$. If
$$
d\geqslant {\rm sep}(G),
$$
where ${\rm sep}(G)$ is the separation
index of the root system of $G$ with
respect to $T$ {\rm(}see Definition~{\rm
\ref{index}} in
Section~{\rm\ref{sepa}}{\rm)},
 then the action of
$G$ on $X_{\lambda_1,\ldots,\lambda_d}$ is
stable {\rm(}and hence ample{\rm)} and the
$G$-stabilizer of a point in general
position in
$X_{\lambda_1,\ldots,\lambda_d}$ is finite.
\end{theorem}

Other results on stability of the
$G$-action on
$X_{\lambda_1,\ldots,\lambda_d}$ are
obtained in Theorem \ref{>>>>} in Section
\ref{stabili}; they are based on some
results from \cite{V}. We show that
stability imposes some constraints on
configuration of the set
$\{\lambda_1,\ldots, \lambda_d\}$ and link
the problem with some monoids that
generalize Littelwood--Richardson
semigroups \cite{zel} whose investigation
during the last decade culminated in
solving several old problems, in
particular, proving Horn's conjecture,
cf.\,survey \cite{fulton}.

 We
apply Theorems \ref{teo1}--\ref{fntn} to
studying primitive $d$-tuples.
Theorem~\ref{teo2}, Definition
\ref{simple}, and Corollary of Theorem
\ref{teo1} immediately imply the following
saturation property.

\begin{theorem}
Let $(\lambda_1,\ldots,\lambda_d)\in {\rm
P}_{\gg}^d$ and $(m_1,\ldots, m_d)\in
\mathbb Z_{>0}^d$. If condition
{\rm\eqref{mff}} holds, then the following
properties are equivalent:
\begin{enumerate}
\item[\rm (i)]  $(\lambda_1, \ldots,
\lambda_d)$ is primitive; \item[\rm (ii)]
$(m_1\lambda_1, \ldots, m_d\lambda_d)$ is
primitive.
\end{enumerate}
\end{theorem}

Utilizing Theorems \ref{teo1} and
\ref{bounds} we prove the following
finiteness theorem about length of
primitive $d$-tuples:

\begin{theorem}\label{fin} Let
$G$ be a simple group. Then for every
primitive $d$-tuple in~${\rm P}_{\gg}^d$,
$$d\leqslant
{\rm sep}(G)+1.$$
\end{theorem}

In view of inequality \eqref{sG} below this
implies

\vskip 2mm

\noindent {\bf Corollary.} {\it Let $G$ be
a simple group. Then $d\leqslant |W|+1$ for
every primitive $d$-tuple in~${\rm
P}_{\gg}^d$.} \vskip 2mm

From Theorem \ref{fin}   we deduce that for
every simple group $G$,
\begin{equation*}\label{pG}
{\rm prim}(G):={\rm sup}\{d\in \mathbb
Z_{>0} \mid \mbox{${\rm P}_{\gg}^d$
contains a primitive element}\}
\end{equation*}
is a natural number not exceeding ${\rm
sep}(G)+1$.
 Discussion in Section 2 and
 the last Corollary
 imply that
$2\leqslant {\rm prim}(G)\leqslant
|W|+1$.

\begin{example} By Example~\ref{e1} we have
${\rm prim}({\bf SL}_2)=3$.
Theorem~\ref{www} below implies that if $G$
is  respectively of type ${\sf A}_l$, ${\sf
B}_l$, ${\sf C}_l$, ${\sf D}_l$, ${\sf
E}_6$, and ${\sf E}_7$, then ${\rm
prim}(G)\geqslant l+2$, $3$, $3$, $3$, $4$,
and  $3$.
\end{example}

For $G={\bf SL}_n$, we can apply to our
problem the representation theory of
quivers. This leads to a characterization
of primitive $d$-tuples of fundamental
weights in terms of canonical decomposition
of dimension vectors of representations of
some graphs and yields a fast al\-gorithm
for verifying whether such a $d$-tuple is
primitive  or not (see
Theorem~\ref{sl-fund} and discussion in
Section \ref{12}).

From Theorems~\ref{teo1}, \ref{multi}, and
\ref{teo2} we deduce a complete
classification of primitive $d$-tuples of
form~\eqref{mfund}:

\begin{theorem}\label{www} Let $G$ be a simple group.
A $d$-tuple
$$(m_1\varpi_i,\ldots,m_d\varpi_i),\hskip 2mm
\mbox{where $d\geqslant 3$,  $(m_1,\ldots,
m_d)\in \mathbb Z_{>0}^d$},$$ is primitive
if and only if the conditions specified in
Table~$2$ hold.
\end{theorem}

Combining Theorem~\ref{teo1} with the
results of \cite{littelmann}, \cite{MWZ1},
and \cite{MWZ2}, we prove that the
following $3$-tuples are primitive.

\begin{theorem} \label{concr}
Let $(\lambda_1,\lambda_2,\lambda_3)\in{\rm
P}_{\gg}^d$. Put $s_i:={\rm
supp}(\lambda_i)$. Then
$(\lambda_1,\lambda_2,\lambda_3)$ is
primitive  in either of the following
cases:

\vskip 4mm

\begin{center}

{\sc Table 3} \vskip 3mm

\begin{tabular}{c|c|c}
$\mbox{\rm no.}$&$\mbox{\rm type of $G$}$ &
{\rm condition}
\\
\hline \hline
 &
 \\[-10pt]
$\begin{matrix} 1\\
2\\
3\\
4\\
5\\
6\\
7
\end{matrix}$
 &
 ${\sf A}_l$ &
 $\begin{matrix}
s_1=\{1\}\\[0pt]
 |s_1|=|s_2|=1\\[0pt]
 |s_1|=|s_2|=1; |s_3|=2\\[0pt]
 |s_1|=1, |s_2|=2; |s_3|=3\\[0pt]
|s_1|=1; |s_2|=2; |s_3|=4\\[0pt]
s_1=\{2\}; |s_2|=2; |s_3|\geqslant 2\\[0pt]
\hskip 6mm |s_1|=1; s_2=\{i, i+1\} \mbox{
\rm or } \{1, j\}, i<l, j\neq 1;
|s_3|\geqslant 2\hskip 6mm {}
\\[-12pt]
&
 \end{matrix}$
 \\[10pt]
\hline
 &\\[-10pt]
$\begin{matrix} 8\\
9\\[-10pt]
\
\end{matrix}$
& ${\sf B}_l$ &
 $\begin{matrix}
\hskip 28.2mm s_1=\{1\}; s_2=\{1,\ldots,l\}
\hskip 28.2mm {}
\\[0pt]
s_1=s_2=\{l\}; s_3=\{1,\ldots,l\}\\[-8pt]
&
 \end{matrix}$
 \\[-3pt]
\hline
 &\\[-10pt]
 $\begin{matrix} 10\\
11\\
12\\
13\\
14\\[-8pt]
\ \end{matrix}$ & ${\sf C}_l$ &
 $\begin{matrix}
 s_1=s_3=\{l\}\\[0pt]
 s_1=\{l\}; s_2=\{i\}, i\neq l; s_3=\{j\},
 j\neq l
 \\[0pt]
\hskip 12.1mm s_1=\{l\}; s_2=\{i\}, i\neq
l; s_3=\{j, m\},
 j\neq m \hskip 12.1mm {}
 \\[0pt]
 s_1=\{l\}; s_2=\{1\}; s_3\neq \{l\}
 \\[0pt]
 s_1=\{1\}; s_2=\{i\}, i\neq l;
 s_3\neq \{l\}
 \\[-6pt]
&
 \end{matrix}$
 \\[-5pt]
\hline
 &\\[-10pt]
 $\begin{matrix} 15\\
16\\
17\\
18\\[-10pt]
\
\end{matrix}$ & ${\sf D}_l$ &
 $\begin{matrix}
s_1=\{1\}; s_2=\{1,\ldots, l\}\\[0pt]
\hskip 17.8mm s_1=\{l-1\}; s_2=\{l\};
s_3=\{1,\ldots,
 l\}\hskip 17.8mm {}
 \\[0pt]
s_1=\{l\}; s_2=\{l\}; s_3=\{1,\ldots,
 l\}\\[0pt]
 s_1=\{3\};  s_2=\{l\}; s_3=\{1,\ldots,
 l\}\\[-8pt]
&
 \end{matrix}$
\\[-3pt]
\hline
 &\\[-10pt]
 $\begin{matrix} 19\\[-10pt]
 \
\end{matrix}$
&$\begin{matrix} {\sf E}_6\\[-10pt]
 \
\end{matrix}$
&
 $\begin{matrix}
 s_1=\{1\}; s_2=\{i\}, i\neq 4;
 s_3=\{1,\ldots,6\}\\[-8pt]
&
 \end{matrix}$
 \\[-3pt]
 \hline
 &\\[-10pt]
 $\begin{matrix} 20
\end{matrix}$
& ${\sf E}_7$ &
 $\begin{matrix}
 s_1=\{1\} \mbox{ \rm or } \{2\}
 \mbox{ \rm or } \{7\};
 s_2=\{7\}; s_3=\{1,\ldots,7\}
 \end{matrix}$\\[-10pt]
 &
\end{tabular}

\vskip 2.5mm
\end{center}
\end{theorem}

Finally, Theorem \ref{geom} below explains
geometric meaning of invariant-freeness of
$d$-tuples and establishes a saturation
property for them. Theorem \ref{bound}
shows that invariant-freeness of
$(\lambda_1,\ldots,\lambda_d)\in{\rm
P}_{\gg}^d$ imposes some constraints on
configuration of the set
$\{\lambda_1,\ldots,\lambda_d\}$ and gives
an upper bound of length of invariant-free
$d$-tuples.

\begin{theorem}\label{geom}
Let $(\lambda_1,\ldots,\lambda_d)\in {\rm
P}_{\gg}^d$ and $(m_1,\ldots, m_d)\in
\mathbb Z_{>0}^d$. The following properties
are equivalent:
\begin{enumerate}
\item[\rm(i)]
$(\lambda_1,\ldots,\lambda_d)$ is
invariant-free;
 \item[\rm (ii)]
$(m_1\lambda_1,\ldots, m_d\lambda_d)$ is
invariant-free; \item[\rm(iii)] the closure
of every $G$-orbit in
$X_{\lambda_1,\ldots,\lambda_d}$ contains
$(0,\ldots,0)\in
E_{\lambda_1}\oplus\ldots\oplus
E_{\lambda_d}$; \item[\rm(iv)]
$k[X_{\lambda_1,\ldots,\lambda_d}]^G= k$.
\end{enumerate}
\end{theorem}

\begin{theorem} \label{bound}
Let $(\lambda_1,\ldots,\lambda_d)\in {\rm
P}_{\gg}^d$ be an invariant-free $d$-tuple.
Then
\begin{enumerate}
\item[\rm(i)] $\mathbb
Q_{>0}\lambda_i^*\notin {\rm
cone}(\{\lambda_1,\ldots,\widehat{\lambda_i},
\ldots,\lambda_d\})$ for every $i$;
\item[\rm(ii)] $d\leqslant {\rm rk}(G)$ if
$G$ is a simple group.
\end{enumerate}
\end{theorem}

\end{enavant}


\begin{enavant} {\bf Primitiveness
and open orbits}\label{po}\nopagebreak

Let $(\lambda_1,\ldots,\lambda_d)\in {\rm
P}_{\gg}^d$. In this section we establish a
connection between the primitiveness of
$(\lambda_1,\ldots,\lambda_d)$ and some
properties of the $G$-actions on
$X_{\lambda_1,\ldots, \lambda_d}$ and
$G/P_{\lambda_1}\times\ldots\times
G/P_{\lambda_d}$.

By \cite[
Theorem\,1]{popov-vinberg1} the variety
$X_{\lambda_i}$ is a cone in $E_{\lambda_i}$, i.e.,
is stable with respect to the action of ${\bf G}_m$
on $E_{\lambda_i}$ by scalar multiplications, and
\begin{equation}\label{boun}
X_{\lambda_i}= {\mathcal O}_{\lambda_i}\cup
\{0\}.
\end{equation}
 This ${\bf G}_m$-action
commutes with the $G$-action and determines
a $G$-stable $k$-algebra $\mathbb
Z^{}_{\geqslant 0}$-grading of
$k[X_{\lambda_i}]$:
\begin{equation}\label{gradu}\textstyle
k[X_{\lambda_i}]=\bigoplus_{n_i\in \mathbb
Z_{\geqslant
0}}k[X_{\lambda_i}]_{n_i},\end{equation}
where $k[X_{\lambda_i}]_{n_i}$ is the space
of ${\bf G}_m$-semi-invariants of the
weight $t\mapsto t^{n_i}$. By
\cite[
Theorem\,2]{popov-vinberg1} there is an isomorphism
of $G$-modules
\begin{equation}\label{isomod}
k[X]_{n_i}\simeq E_{n_i\lambda_i^*}.
\end{equation}

The group $G^d\times {\bf G}_m^d$ acts on
$X_{\lambda_1,\ldots,\lambda_d}$ in the
natural way, and $\mathcal
O_{\lambda_1,\ldots,\lambda_d}$ is a ${\bf
G}_m^d$-stable open $G^d$-orbit in
$X_{\lambda_1,\ldots,\lambda_d}$. By
restriction this action entails the action
of $G\times {\bf G}_m^d$, where $G$ is
diagonally embedded in $G^d$. The action of
${\bf G}_m^d$ on
$X_{\lambda_1,\ldots,\lambda_d}$ determines
a $G$-stable $k$-algebra $\mathbb
Z_{\geqslant 0}^d$-grading of
$k[X_{\lambda_1,\ldots,\lambda_d}]$,
\begin{equation}\label{mgradu}\textstyle
k[X_{\lambda_1,\ldots,\lambda_d}]=
\bigoplus_{n_1,\ldots, n_d\in \mathbb
Z_{\geqslant 0}}
k[X_{\lambda_1,\ldots,\lambda_d}]_{(n_1,\ldots,
n_d)},\end{equation} where
$k[X_{\lambda_1,\ldots,\lambda_d}]_{(n_1,\ldots,
n_d)}$ is the space of ${\bf
G}_m^d$-semi-invariants of the weight
$(t_1,\ldots, t_d)\!\mapsto t_1^{n_1}\cdots
t_d^{n_d}$. Since, by \eqref{X},
 $k[X_{\lambda_1,\ldots,\lambda_d}]$ and
$\bigotimes_{i=1}^{d}k[X_{\lambda_i}]$ are
$G$-isomorphic $k$-algebras,
 \eqref{gradu}, \eqref{isomod},
 and \eqref{mgradu}
yield that, for every $(n_1,\ldots, n_d)\in
\mathbb Z_{\geqslant 0}^d$, there is an
isomorphism of $G$-modules
\begin{equation}\label{isomo}\textstyle
k[X_{\lambda_1,\ldots,\lambda_d}]_{(n{}^{}_1,\ldots,
n_d)}\simeq\bigotimes_{1\leqslant
i\leqslant d}k[X_{\lambda_i}]_{n_i}\simeq
\bigotimes_{1\leqslant i\leqslant
d}E_{n{}^{}_i\lambda_i^*}.
\end{equation}

Consider now the categorical quotient
\eqref{cat} for $H=G$ and
$Y=X_{\lambda_1,\ldots,\lambda_d}$ and
denote
 $\pi^{}_{H, Y}$ by
$\pi_{\lambda_1,\ldots,\lambda_d}$. The
field of fractions of
$k[X_{\lambda_1,\ldots,\lambda_d}]^G$ is
then
$\pi_{\lambda_1,\ldots,\lambda_d}^*(k(\XC))$.
Since the action ${\bf G}_m^d$ on
$X_{\lambda_1,\ldots,\lambda_d}$ commutes
with that of $G$, it descends to $\XC$. The
corresponding action of ${\bf G}_m^d$ on
the $k$-algebra $k[\XC]$ determines its
$\mathbb Z_{\geqslant 0}^d$-grading
\begin{gather}\label{mgraduXC}
\begin{gathered}\textstyle
k[\XC]= \bigoplus_{n_1,\ldots, n_d\in
\mathbb Z_{\geqslant 0}}
k[\XC]_{(n_1,\ldots, n_d)},\quad
\mbox{where}\\
\pi_{\lambda_1,\ldots,\lambda_d}^*\bigl(k[\XC]_{(n_1,\ldots,
n_d)}\bigr)=k[X_{\lambda_1,\ldots,\lambda_d}]^G\cap
k[X_{\lambda_1,\ldots,\lambda_d}]_{(n_1,\ldots,
n_d)}.\end{gathered} \end{gather} From
this, \eqref{mgradu}, and \eqref{isomo} we
deduce that
\begin{equation}\label{XGG}\textstyle
k[\XC]_{(n_1,\ldots,n_d)}\simeq
\bigl(\bigotimes_{1\leqslant i\leqslant
d}E_{n{}^{}_i\lambda_i^*}\bigr)\hskip
-.1mm^{G}.
\end{equation}

\begin{lemma}\label{crsimpl}
The following properties are equivalent:
\begin{enumerate}
\item[\rm(i)] $k(\XC)^{{\bf G}_m^d}=k$;
\item[\rm(ii)] there is an open ${\bf
G}_m^d$-orbit in $\XC$; \item[\rm(iii)]
$(\lambda_1,\ldots,\lambda_d)$ is
primitive.
\end{enumerate}
\end{lemma}
\begin{proof}
(i)$\Leftrightarrow$(ii) follows from
Rosenlicht's theorem \cite{rosenlicht},
cf., e.g.,\,\cite[Cor.\,of 
Theo\-rem 2.3]{popov-vinberg2}.

Assume that $(\lambda_1,\ldots,\lambda_d)$
is not primitive. Then $\dim
\bigl(\bigotimes_{1\leqslant i\leqslant
d}E_{n_i\lambda_i} \bigr)\hskip
.0mm^G\geqslant 2$ for some $(n_1,\ldots,
n_d)\in\mathbb Z_{\geqslant 0}^d$. Since
for all $(\mu_1,\ldots, \mu_d)\in{\rm
P}_{++}^d$ and $(m_1,\ldots, m_d)\in
\mathbb Z_{\geqslant 0}^d$, we have
\begin{gather}\label{**}\begin{gathered}\textstyle
\dim \bigl( \bigotimes_{1\leqslant
i\leqslant d}E_{m_i\mu_i} \bigr)\hskip
.0mm^G= \dim\bigl(\bigotimes_{1\leqslant
i\leqslant d}E_{m^{}_i\mu_i^*} \bigr)\hskip
.0mm^G,
\end{gathered}
\end{gather}
it then follows from \eqref{XGG} that $\dim
(k[\XC]_{(n_1,\ldots,n_d)})\geqslant 2$.
This means that the algebra $k[\XC]$
contains two nonproportional ${\bf
G}_m^d$-semi-invariant functions $f_1$ and
$f_2$ of the same weight. Hence $f_1/f_2
\in k(\XC)^{{\bf G}_m^d}$, $f_1/f_2\notin
k$. This proves (i)$\Rightarrow$(iii).

Conversely, let $f\in k(\XC)^{{\bf
G}_m^d}$, $f\notin k$. Since $\XC$ is an
affine variety, $k(\XC)$ is the field of
fractions of $k[\XC]$. As ${\bf
 G}_m^d$ is a connected solvab\-le group,
 by
\cite[
Theorem\,3.3]{popov-vinberg2}  this implies that in
$k[\XC]$ there are two ${\bf
 G}_m^d$-semi-invariant elements
of the same ${\bf
 G}_m^d$-weight, say, $f_1, f_2\in
 k[\XC]_{(n_1,\ldots,n_d)}$,
such that $f=f_1/f_2$. Since $f_1$ and
$f_2$ are nonproportional, $\dim
(k[\XC]_{(n_1,\ldots,n_d)})\geqslant 2$. By
\eqref{XGG} and \eqref{**} this yields
$\dim\bigl( \bigotimes_{i=1}^d
E_{n{}^{}_i\lambda_i}\bigr)\hskip
-.1mm^G\geqslant 2$. Hence
$(\lambda_1,\ldots,\lambda_d)$ is not
primitive. This proves
(iii)$\Rightarrow$(i).
 \quad $\square$
\renewcommand{\qed}{}
\end{proof}
\begin{remark} \label{toric} By
\cite[
Theorem\,3]{popov-vinberg2} every $X_{\lambda_i}$ is
a normal variety. From \eqref{X} we then conclude
that $X_{\lambda_1,\ldots, \lambda_d}$ is normal as
well. Hence property (ii) in Lemma \ref{crsimpl}
means that $\XC$ is a toric ${\bf G}_m^d$-variety.
\end{remark}

\begin{lemma}\label{open} The following
properties are equivalent:
\begin{enumerate}
\item[\rm(i)]
$k(X_{\lambda_1,\ldots,\lambda_d})^{G\times
{\bf G}_m^d}=k$; \item[\rm(ii)] there is an
open $G\times{\bf G}_m^d$-orbit in
$X_{\lambda_1,\ldots,\lambda_d}$;
\item[\rm(iii)] there is an open $G$-orbit
in $ G/P_{\lambda_1} \times\ldots\times
G/P_{\lambda_d}$.
\end{enumerate}
\end{lemma}
\begin{proof} The aforementioned
Rosenlicht's theorem yields the
equivalencies (i)$\Leftrightarrow$(ii) and
\begin
{gather}\label{f2}
\begin{gathered}
k\bigl(G/P_{\lambda_1}\!\times\ldots\times\!
G/P_{\lambda_d}\bigr)\hskip -.3mm^G \hskip
-1mm=\!k \hskip -.9mm\iff\hskip -.9mm
\mbox{$\displaystyle G/P_{\lambda_1}
\!\times\ldots\times \!G/P_{\lambda_d}$
contains an open $G$-orbit.}
\end{gathered}
\end{gather}

The natural projection
$\rho_{\lambda_i}\colon \mathcal
O_{\lambda_i}\rightarrow G/P_{\lambda_i}$
is a rational quotient for the ${\bf
G}_m$-action on $\mathcal O_{\lambda_i}$,
cf.\;\cite[2.4]{popov-vinberg2}. Hence the
$G$-equivariant morphism
$\rho_{\lambda_1}\times\ldots\times
\rho_{\lambda_d}\colon \mathcal
O_{\lambda_1,\ldots,\lambda_d} \rightarrow
G/P_{\lambda_1}\times\ldots\times
G/P_{\lambda_d}$ is a rational quotient for
the ${\bf G}_m^d$-action on $\mathcal
O_{\lambda_1,\ldots, \lambda_d}$. Therefore
it induces an isomorphism of invariant
fields
\begin{equation}\label{f1}
  k\bigl(
  G/P_{\lambda_1}
  \times\ldots\times
G/P_{\lambda_d} \bigr)\hskip -.1mm^G
 \overset{\simeq}{\longrightarrow}
 k(\mathcal O_{\lambda_1,\ldots,\lambda_d})^{G\times {\bf
 G}_m^d}.
\end{equation}

But $k(\mathcal
O)=k(X_{\lambda_1,\ldots,\lambda_d})$ since
$\mathcal O$ is open in
$X_{\lambda_1,\ldots,\lambda_d}$. This,
\eqref{f2}, and \eqref{f1} now imply
(i)$\Leftrightarrow$(iii). \quad $\square$
\renewcommand{\qed}{}
\end{proof}

\begin{lemma}\label{semigr}
Let $(\lambda_1,\ldots, \lambda_d)\in {\rm
P}_{\gg}^d$, $s\in \mathbb Z_{> 0}$, and
$${\mathcal M}_d(s):=\{
(n_1,\ldots, n_d)\in\mathbb Z_{\geqslant
0}^d\mid c_{n_1\lambda_1,\ldots,
n_d\lambda_d}^0\geqslant s\}.$$ Then
\begin{equation*}\label{ideal}
{\mathcal M}_d(1)+
{\mathcal M}_d(s)\subseteq {\mathcal
M}_d(s).
\end{equation*}
\end{lemma}
\begin{proof} By \eqref{XGG} and \eqref{**}
\begin{equation}\label{mmm}
(n_1,\ldots, n_d)\in {\mathcal
M}_d(s)\hskip 2mm \iff\hskip 2mm
\dim\bigl(k[\XC]_{(n_1,\ldots,n_d)}\bigr)\geqslant
s.
\end{equation}
Let $\alpha\in{\mathcal M}_d(1)$ and
$\beta\in {\mathcal M}_d(s)$. Pick a
nonzero function $f\in k[\XC]_{\alpha}$ and
linear independent functions $h_1,\ldots,
h_s\in k[\XC]_{\beta}$: by \eqref{mmm} this
is possible. Then the functions
$fh_1,\ldots, fh_d$ are li\-nearly
independent since $k[\XC]$ is an integral
domain.  They lie $k[\XC]_{\alpha+\beta}$
since \eqref{mgraduXC} is a grading. So
$\dim\bigl(k[\XC]_{\alpha+\beta}\bigr)\geqslant
s$, whence $\alpha+\beta\in {\mathcal
M}_d(s)$ by \eqref{mmm}. \quad $\square$
\renewcommand{\qed}{}
\end{proof}
\end{enavant}


\begin{enavant}{\bf  Proof of Theorem~{\bf
\ref{teo1}}}

If the assumption of (i) holds,
Lemma~\ref{open} implies that
\begin{equation*}
k(X_{\lambda_1,\ldots,\lambda_d})^{G
\times{\bf G}_m^d}=\bigl(k(X_{\lambda_1,
\ldots,\lambda_d})^G\bigr)\hskip
-.1mm^{{\bf G}_m^d}=k.\end{equation*} Since
\begin{equation}\label{subset}
\pi_{\lambda_1,\ldots, \lambda_d}^*\colon
k(\XC)\hookrightarrow
k(X_{\lambda_1,\ldots,\lambda_d})^G,
\end{equation}
this yields $k(\XC)\hskip -.1mm^{{\bf
G}_m^d}=k$; whence
$(\lambda_1,\ldots,\lambda_d)$ is primitive
by Lemma~\ref{crsimpl}. This proves~(i).

If the assumption of (ii) holds, let
\begin{equation*}\varrho\colon
X_{\lambda_1,\ldots,\lambda_d}
\dashrightarrow \XG\end{equation*} be a
rational quotient for the action of $G$ on
$X_{\lambda_1,\ldots,\lambda_d}$. i.e.,
$\XG$ is an irreducible variety and
$\varrho$  a dominant rational map such
that $\varrho^*\bigl(k(\XG)\bigr)=
k(X_{\lambda_1,\ldots,\lambda_d})^G$, cf.,
e.g.,\,\cite[2.4]{popov-vinberg2}. By
\cite[Prop.\,2.6]{popov-vinberg2} the
action of ${\bf G}_m^d$ on
$X_{\lambda_1,\ldots,\lambda_d}$
 induces a
rational ${\bf G}_m^d$-action on $\XG$ such
that $\varrho$ becomes ${\bf
G}_m^d$-equivariant.  By \cite[Cor.\,of
Theorem\,1.1]{popov-vinberg2} replacing $\XG$ with a
birationally isomorphic variety, we may (and shall)
assume that the rational action of ${\bf G}_m^d$ on
$\XG$ is regular (morphic).

Embedding \eqref{subset} induces a dominant
rational ${\bf G}_m^d$-equivari\-ant map
$\tau\colon \XG\dashrightarrow \XC$ such
that we obtain the following commutative
diagram:
\begin{equation}\label{ratio}
\begin{matrix}
\xymatrix{&X_{\lambda_1,\ldots,
\lambda_d}\ar@{-->}[dl]_{\varrho}
\ar[dr]^{\pi_{\lambda_1,\ldots, \lambda_d}}&\\
\XG\ar@{-->}[rr]^{\tau}&&\XC}
\end{matrix}\quad.
\end{equation}

\vskip 3mm

Since the action of $G$ on
$X_{\lambda_1,\ldots, \lambda_d}$ is ample,
we have
\begin{equation}\label{=}
\dim(\XG)=\dim(\XC).
\end{equation}
Since $(\lambda_1,\ldots,\lambda_d)$ is
primitive, Lemma~\ref{crsimpl} yields that
$\XC$ contains an open ${\bf G}_m^d$-orbit.
From this, \eqref{ratio}, and \eqref{=} it
then follows that $\XG$ contains an open
${\bf G}_m^d$-orbit. Hence
$\bigl(k(\XG)\bigr)\hskip -.1mm^{{\bf
G}_m^d}=k$, i.e.,
$\bigl(k(X_{\lambda_1,\ldots,\lambda_d})^G
\bigr)\hskip -.1mm^{{\bf
G}_m^d}=k(X_{\lambda_1,\ldots,\lambda_d})\hskip
-.1mm^{G\times{\bf G}_m^d}=k$.
Lemma~\ref{open} then implies that $
G/P_{\lambda_1} \times\ldots\times
G/P_{\lambda_d}$ contains an open
$G$-orbit. This proves (ii).\quad $\square$
\end{enavant}


\begin{enavant} {\bf Proof of Theorem
\ref{bounds}}

\nopagebreak Since dimension of a multiple
flag variety $G/P_{\lambda_1}\times\ldots
\times G/P_{\lambda_d}$ containing an open
$G$-orbit does not exceed $\dim(G)$, we
have $$\dim (G)\geqslant \sum_{i=1}^d
\bigl(\dim (G)-\dim
(P_{\lambda_i})\bigr)\geqslant d\bigl(\dim
(G)-\dim (P)\bigr),$$ where $P$ is a
parabolic subgroup of $G$ of maximal
dimension. Since $\dim (G)=\dim (L)+2\dim
(P_u)$ and $\dim (P)=\dim (L)+\dim (P_u)$,
where $L$ and $P_u$ are respectively a Levi
subgroup and the unipotent radical of $P$,
this yields
\begin{equation}\label{estimate1}
d\leqslant \frac{2\dim G}{\dim G-\dim L}.
\end{equation}

Let $L_i$ be a Levi subgroup of
$P_i:=P_{\varpi_i}$. The equality $\dim
(L_i)=2\dim (P_{i})-\dim(G)$ implies that
$\dim (L_i)$ and $\dim (P_i)$, as functions
in $i$, attain their absolute maximums at
the same values of $i$; let $M$ be the set
of these values. Then the group $P$ is
conjugate to $P_{i_0}$ for some $i_0\in M$.

Since $L_i$ is a reductive group of rank
${\rm rk}(G)$ and the Dynkin diagram of its
commutator group is obtained from that of
$G$ by removing the $i$th node, finding all
the $\dim (L_i)$'s and then the set $M$ is
a matter of some clear calculations. We
skip them (see some details in
\cite[Sect.\,5--13]{popov2}). The results
are collected in Table~4 below.

\vskip 3mm

\begin{center}
{\sc Table 4}

\vskip 2mm
\begin{tabular}{c||c|c|c|c|c|c|c|c|c}
\text{ {\rm type of} $G$} & ${\sf A}_l$,
\mbox{\fontsize{8pt}{5mm}\selectfont
$l\geqslant 1$}& ${\sf B}_l,$
\mbox{\fontsize{8pt}{5mm}\selectfont
$l\geqslant 3$}& ${\sf C}_l,$
\mbox{\fontsize{8pt}{5mm}\selectfont
$l\geqslant 2$}& ${\sf D}_l$,
\mbox{\fontsize{8pt}{5mm}\selectfont
$l\geqslant 4$}& ${\sf E}_6$ & ${\sf E}_7$&
${\sf E}_8$& ${\sf F}_4$& ${\sf G}_2$
\\[2pt]
 \hline
&&&&&&&&&\\[-11pt]
$M$ &$1$, $l$&$1$&$1$&$1$&$1$, $6$&$7$&$8$&$1$, $4$&$1$, $2$\\[1pt]
\hline
&&&&&&&&&\\[-9pt]\\[-15pt]
\mbox{\fontsize{13pt}{5mm}\selectfont
$\frac{2\dim G}{\dim G-\dim L}$}
&$l+2$&\mbox{\fontsize{13pt}{5mm}\selectfont
$\frac{2l^2+l}{2l-1}$}&\mbox{\fontsize{13pt}{5mm}\selectfont
$\frac{2l^2+1}{2l-1}$}&\mbox{\fontsize{13pt}{5mm}\selectfont
$\frac{2l^2-l}{2l-2}$}&\mbox{\fontsize{13pt}{5mm}
\selectfont
$\frac{39}{8}$}&\mbox{\fontsize{13pt}{5mm}
\selectfont
$\frac{133}{27}$}&\mbox{\fontsize{13pt}{5mm}\selectfont
$\frac{248}{57}$}&\mbox{\fontsize{13pt}{5mm}\selectfont
$\frac{52}{15}$}&\mbox{\fontsize{13pt}{5mm}\selectfont
$\frac{14}{5}$}
\end{tabular}\\
\vskip 3mm
\end{center}
The claim now immediately follows from
\eqref{estimate1} and Table 4.\quad
$\square$

\end{enavant}


\begin{enavant}{\bf  Proof of Theorem~{\bf
\ref{teo2}}}

We put, for brevity,
$$X:=X_{\lambda_1,\ldots,\lambda_d},\hskip
2mm \mathcal O:=\mathcal
O_{\lambda_1,\ldots,\lambda_d},\hskip 2mm
v:=v_{\lambda_1}\times\ldots\times
v_{\lambda_d}\in\mathcal O.$$ Clearly,
$\dim (\overline{\mathcal
O}_{\lambda})\geqslant 2$ for every
$\lambda \in  {\rm P}_\gg$, hence by
\eqref{X} and \eqref{boun}
\begin{equation}\label{bbb}
{\rm codim}^{}_{X}(X\setminus \mathcal
O)\geqslant 2.
\end{equation}
Recall from Remark \ref{toric} that $X$ is
normal. From \eqref{bbb} we conclude that
\begin{equation}\label{ClXO}
{\rm Cl}(X)\simeq {\rm Cl}(\mathcal O).
\end{equation}

Since
$G^d_v=G_{v_{\lambda_1}}\times\ldots\times
G_{v_{\lambda_d}}$ and $G^d$ is a connected
simply connected semisimple group, we
deduce from \cite[Prop.\,1; Cor. of
Theorem\,4]{popov1} that
\begin{equation}\label{Picard}
{\rm Cl}(\mathcal O)\simeq {\rm
Cl}(G^d/G^d_v)\simeq \mathcal
X(G^d_v)\simeq \mathcal
X(G_{v_{\lambda_1}})\oplus\ldots\oplus
\mathcal X(G_{v_{\lambda_d}}).
\end{equation}
On the other hand, \eqref{mff} and
\cite[\S1, no.\,5]{popov-vinberg1} imply
that
\begin{equation}\label{chara}
\mathcal X(G_{v_{\lambda_i}})\simeq \mathbb
Z/m_i.
\end{equation}
From \eqref{ClXO}, \eqref{Picard}, and
\eqref{chara} we obtain that
\begin{equation}\label{classx}\textstyle
{\rm Cl}(X)\simeq \bigoplus_{i=1}^d \mathbb
 Z/m_i.
\end{equation}

Further, since $G$ is semisimple, we have
\begin{equation}\label{xi}
\mathcal X(G)=\{0\};\end{equation} whence
every invertible element of $k[G]$
 is constant, see \cite{rosenlicht0}.
  Therefore the same holds for
 $k[\mathcal O]$ as well.
As $\mathcal O$ is open in $X$, this yields
that every in\-ver\-tible element of $k[X]$
is constant.

Take now  a nonconstant function $f\in
k(X)^G$. Then $(f)\neq 0$. For, otherwise,
the normality of $X$ would imply (see,
e.g.,\,\cite[
Theorem\,38]{M}) that $f$ is invertible element of
$k[X]$, hence a constant, a contradiction.

Since, by \eqref{classx},  ${\rm Cl}(X)$ is
a finite group, there is $n\in\mathbb
Z^{}_{> 0}$ such that both divisors
$n(f)_0$ and $n(f)_\infty$ are principal,
i.e.,
\begin{equation}\label{h1h2}
\mbox{$n(f)_0=(h_1)$ and
$n(f)_\infty=(h_2)$ for some $h_1, h_2\in
k(X)$.}
\end{equation} As $n(f)_0\geqslant 0$ and
$n(f)_\infty\geqslant 0$, the normality of
$X$ and \eqref{h1h2} imply
 that $h_1,
h_2\in k[X]$ (see,
e.g.,\,\cite[
Theorem\,38]{M}). Further, since $f$ is
$G$-invariant, the supports of $(f)_0$ and
$(f)_\infty$ are $G$-stable subsets of $X$. By
\cite[Theorem 3.1]{popov-vinberg2} this and
\eqref{h1h2} imply that $h_1$ and $h_2$ are
$G$-semi-invariants. Hence by \eqref{xi}
\begin{equation}\label{invv}
h_1, h_2\in k[X]^G.
\end{equation}
On the other
hand, $(f^nh_2/h_1)=0$ by \eqref{h1h2},
hence $f^nh_2/h_1$ is a constant. By
\eqref{invv}
 this means that $f$ is algebraic over
 the field of fractions of
 $k[X]^G$. Hence, by
 Definition~\ref{ample},
 the action of $G$ on $X$ is ample. This
 completes the proof.\quad $\square$

 \end{enavant}


\begin{enavant}{\bf Separation index
of irreducible
root
system}\label{sepa}

Let $R$ be a
root system in a rational vector
space $L$
(we assume that $L$ is the linear span of $R$)
and let
$W(R)$ be the Weyl group of $R$. For any linear
function $l\in L^*$, put
\begin{equation}\label{+-}
l^+:=\{x\in L\mid l(x)\geqslant 0\},\hskip
2mm l^0:=\{x\in L\mid l(x)=0\}, \hskip 2mm
l^-:=\{x\in L\mid l(x)\leqslant 0\}.
\end{equation}
Given a subset $S$ of $L$, denote by
$\overline S$ the closure of $S$ in $L$.

\begin{lemma}\label{oopen}
Let $R$ be an irreducible root system. Then for every
nonzero linear function $l\in L^*$, there is a Weyl
chamber $C\subset L$ of $R$ such that
\begin{equation*}\label{C}
\overline C\subset l^+ \quad\mbox{and}\quad
\overline C \cap l^0=\{0\}.
\end{equation*}
\end{lemma}

\begin{proof}
First, we prove that $R\cap l^+$ contains a
basis of $R$. If $R\cap l^0=\varnothing$,
this is proved, e.g., in \cite[\S8, Prop.
4]{serre}. In general case, fix  a choice
of Euclidean structure on $L^*$ and let $S$
be a ball in $L^*$ with the center at $l$.
We identify in the natural way every
$\alpha\in R$ with a linear function on
$L^*$. Taking $S$ small enough, we may (and
shall) assume that every $\alpha\in
R\setminus l^0$ has no zeros on $S$. On the
other hand, since $R$ is finite, $S$ does
not lie in the union of hyperplanes defined
by vanishing of the roots from $R\cap l^0$.
Hence there is an element $s\in S$ such
that
\begin{equation}\label{hf}
R\cap s^0=\varnothing\quad \mbox{and} \quad
R\cap l^+\supseteq R\cap s^+.\end{equation}
According to the aforesaid, the equality in
\eqref{hf} implies that $R\cap s^+$
contains a basis of $R$. Then the inclusion
in \eqref{hf} yields the claim.

Let now $\beta_1,\ldots,\beta_r$ be a basis of $R$
contained in $R\cap l^+$. Then
\begin{equation}\label{cone}
\sum_i\mathbb Q_{>0}\beta_i\subset l^+
\setminus l^0.
\end{equation}
Let $\pi_1,\ldots, \pi_r\in L$ be the basis of $L$
dual to $\beta_1^\vee,\ldots,\beta_r^\vee$ (i.e.,
$\pi_1,\ldots, \pi_r\in L$ are the fundamental
weights corresponding to $\beta_1,\ldots,\beta_r$).
Then
\begin{equation}\label{pb}
\pi_i=\sum_j c_{ij}\beta_j,
\end{equation}
where $c_{ij}$ are the elements of inverse
Cartan matrix of $R$. Since $R$ is
irreducible,
\begin{equation}\label{ce}
c_{ij}\in \mathbb Q_{>0}\quad\mbox{for all
} i \mbox{ and } j,
\end{equation}
see, e.g.,\,\cite{OV}. Consider now the
Weyl chamber $C:=\sum_i\mathbb
Q_{>0}\pi_i$. Since $\overline
C:=\sum_i\mathbb Q_{\geqslant 0}\pi_i$, it
follows from \eqref{pb} and \eqref{ce},
that $\overline C\setminus
\{0\}\subset\sum_i\mathbb Q_{>0}\beta_i$.
Now the claim follows from \eqref{cone}.
 \quad
$\square$
\renewcommand{\qed}{}
\end{proof}

\noindent{\bf Corollary.}\label{cor} {\it
Let $R$ be an irreducible root system. Then
there is a sequence $C_1,\ldots, C_n\subset
L$ of the Weyl chambers of $R$ satisfying
the following property:}
\begin{gather}\label{prop}
\begin{gathered} \mbox{\it for every
nonzero linear function $l\in L^*$, there
is a natural $i\in [1,
n]$ such that}\\[-3pt]
\overline {C_i}\subset l^+ \quad\mbox{\it
and}\quad \overline {C_i} \cap l^0=\{0\}.
\end{gathered}
\end{gather}
\begin{proof} Let $C_1,\ldots, C_n$ be
 a sequence
 of all Weyl chambers of $R$.
 By Lemma~\ref{oopen} it
satisfies property \eqref{prop}.
 \quad
$\square$
\renewcommand{\qed}{}
\end{proof}

\begin{definition}\label{index} Let $R$ be an
irreducible
root system in a
rational vector space $L$. The {\it separation index}
${\rm sep}(R)$ of $R$ is the minimal length of
sequences $C_1,\ldots, C_n$ of
 Weyl chambers
 of $R$ satisfying
property \eqref{prop}.
\end{definition}

\begin{lemma}\label{ner}
The following inequalities hold:
\begin{equation}\label{sG}
{\rm rk}(R)+1\leqslant {\rm
sep}(R)\leqslant |W(R)|.
\end{equation}
\end{lemma}
\begin{proof}
Let $C_1,\ldots, C_{{\rm sep}(R)}$ be a
sequence of Weyl chambers of $R$ satisfying
property \eqref{prop}. For every $i$, fix a
choice of point $x_i\in C_i$. Arguing on
the the contrary, assume that  ${\rm
sep}(R)\leqslant {\rm rk}(R)$. Then there
is a nonzero linear function $l\in L^*$
such that $x_i\in l^-$ for all ~$i$. This
contradicts property \eqref{prop}. Thus the
left inequality in \eqref{sG} is proved.
The right one follows from the fact that
$|W(R)|$ is equal to the cardinality of
set of all Weyl chambers of $R$.
 \quad
$\square$
\renewcommand{\qed}{}
\end{proof}

The example below shows that all equalities
and inequalities in \eqref{sG} are attained
for suitable~$R$'s.
\begin{example} Clearly,
${\rm sep}({\sf A}_1)=2$, and it is not
difficult to  verify that
$$
{\rm sep}({\sf A_2})=6,\quad {\rm sep}({\sf
B_2})=4, \quad \mbox{and}\quad {\rm
sep}({\sf G_2})=3
$$
(since ${\rm sep}(R)$ depends only on the
type of $R$, the meaning of notation is
clear).
\end{example}
\end{enavant}

\begin{remark}\label{refl}
The notion of separation index can be defined in a
more general setting.

Namely, let $\mathcal M$ be a finite set of nonempty
subsets of a finite dimensional real vector space
$L$.
\begin{definition} We call
a subset $\mathcal S$ of $\mathcal M$ {\it
separating} for $\mathcal M$ if for every nonzero
linear function $l\in L^*$ there exists a set $M
\in\mathcal S$ such that $l$ is strictly positive at
every nonzero point of $M$. If there exists a
separating set for $\mathcal M$, we say that the {\it
separation property} holds for $\mathcal M$ and call
the minimum ${\rm sep}(\mathcal M)$ of cardinalities
of separating sets for $\mathcal M$ the {\it
separation index} of~$\mathcal M$.
\end{definition}

Arguing as in the proof of Lemma~\ref{ner}, we obtain
$\dim(L)+1\leqslant {\rm sep}(\mathcal M)\leqslant
|\mathcal M|$.

\begin{example} Assume that
\begin{enumerate}
\item[\rm (a)] $L=\bigcup_{M\in \mathcal M} M$;
\item[\rm(b)] there exists a Euclidean inner product
$\langle\ {,}\ \rangle$ on $L$ such that for every
set $M\in \mathcal M$, the angle between every two
nonzero vectors of $M$ is acute.
\end{enumerate} Then the {\it separation
property} holds for $\mathcal M$. Indeed, let $l\in
L^*$, $l\neq 0$. Identify $L^*$ with $L$ by means of
$\langle\ {,}\ \rangle$. Then (a) implies that $l$
lies in some $M_0\in \mathcal M$, and (b) implies
that $l$ is strictly positive at every nonzero point
of $M_0$.
\end{example}

\begin{example}
Let $K\subset {\bf GL}(L)$ be an irreducible finite
reflection group and let $\mathcal M$ be the set of
closures of its chambers in $L$. Then (a) holds. Let
$\langle\ {,}\ \rangle$ be a $K$-invariant Euclidean
inner product on $L$. Then the irreducibility of $K$
implies that  (b) holds as well. Hence, in this case,
the separation property holds for $\mathcal M$.
\begin{definition}\label{sepref}  In this case we call
${\rm sep}(\mathcal M)$  the {\it separation index
of} $K$ and denote it by ${\rm sep}(K)$.
\end{definition}
\end{example}

Definitions \ref{index} and \ref{sepref} imply that
if $K$ is crystallographic, i.e., $K=W(R)$ for an
irreducible root system $R$, then ${\rm sep}(R)={\rm
sep}(K)$. The next example illustrates the
noncrystallographic case.
\begin{example} It is not difficult to verify that
${\rm sep}(I_2(p))=3$ for $p\geqslant 7$ and ${\rm
sep}(I_2(5))=4$.
\end{example}

It would be interesting to calculate ${\rm sep}(K)$ for
every irreducible finite reflection group $K$ and, in
particular, to find ${\rm sep}(\Delta)$ for every
irreducible root system $\Delta$.\footnote{{\it Added in
proof}. Recently in V. Zhgoon, D. Mironov, {\it
Separating systems of Weyl chambers}, Math. Notes, to
appear, the following upper bounds have been obtained:
${\rm sep}({\sf A}_l)\leqslant 2l!+2$, ${\rm sep}({\sf
B}_l)={\rm sep}({\sf C}_l)\leqslant 2^{l+1}-2$, ${\rm
sep}({\sf D}_l)\leqslant 2^{l-1}l!/(l-1)+2$, ${\rm
sep}({\sf F}_4)\leqslant 30$, ${\rm sep}({\sf
E}_6)\leqslant 242$, ${\rm sep}({\sf E}_7)\leqslant
4610$, ${\rm sep}({\sf E}_8)\leqslant 9222$, ${\rm
sep}({\sf H}_3)\leqslant 14$, ${\rm sep}({\sf
H}_4)\leqslant 30$.}
\end{remark}


\begin{enavant} {\bf Proof of Theorem
\ref{fntn}}

The proof is based on the following lemma.
\begin{lemma} \label{closedness} Let $V=V_1\oplus
\ldots\oplus V_n$, where $V_1,\ldots, V_n$
are finite dimensional $G$-modules, let
$v_i\in V_i$ be a $T$-weight vector of a
weight $\mu_i\in \mathcal X(T)$, and let
$v:=v_1+\ldots +v_n\in V$. Then the
following properties  are equivalent:
\begin{enumerate}
\item $G\cdot v$ is closed; \item $T\cdot
v$ is closed; \item $0\in {\rm
int}\bigl({\rm conv}(\{\mu_1,\ldots,
\mu_n\})\bigr)$.
\end{enumerate}
\end{lemma}

\begin{proof} This  is proved in
\cite[
Theorem\,1]{popov3}. \quad $\square$
\renewcommand{\qed}{}
\end{proof}

Passing to the proof of Theorem \ref{fntn},
we first  establish the existence of
elements $w_1,\!\ldots,\! w_d$ of $W$ such
that
\begin{enumerate}
\item[\rm(i)] $\dim \bigl({\rm
conv}(\{w_1\cdot\lambda_1,\ldots,
w_d\cdot\lambda_d\})\bigr)=r (={\rm
rk}(G))$; \item[\rm(ii)] $0\in {\rm
int}\bigl({\rm
conv}(\{w_1\cdot\lambda_1,\ldots,
w_d\cdot\lambda_d\})\bigr)$.
\end{enumerate}

Let $R\subset L:=\mathcal X(T)_{\mathbb Q}$
be the root system of $G$ with respect to
$T$ and let $C_1,\ldots,C_{{\rm sep}(G)}$
be a sequence of Weyl chambers of $R$
satisfying property \eqref{prop}. For every
$i\leqslant {\rm sep}(G)$, let $w_i$ be the
(unique) element of $W$ such that $w_i\cdot
\lambda_i\in \overline{C_i}$. For every
$i\geqslant {\rm sep}(G)+1$, put $w_i=e$.
If (i) or (ii) fails, then ${\rm
conv}(\{w_1\cdot\lambda_1,\ldots,
w_d\cdot\lambda_d\})\subset l^-$ for some
linear function $l\in L^*$. But the choice
of $C_1,\ldots, C_{{\rm sep}(G)}$ implies
that there is $i\leqslant {\rm sep}(G)$
such that $\overline{C_i}\subset l^+$ and
$\overline{C_i}\cap l^0=\{0\}$. Since
$\lambda_i\neq 0$, we have
$w_i\cdot\lambda_i\in
\overline{C_i}\setminus\{0\}$, hence
$w_i\cdot\lambda_i\in l^+\setminus l^0$.
Therefore $w_i\cdot\lambda_i\notin l^-$, a
contradiction. Thus (i) and (ii) hold, and
the existence of desired
 $w_i$'s
 is proved.

 Consider now
the point
\begin{equation}\label{point}
v:=\overset{.}w_1\cdot v_{\lambda_1}+\ldots
+ \overset{.}w_d\cdot v_{\lambda_d} \in
X_{\lambda_1,\ldots, \lambda_d}\subseteq
E_{\lambda_1}\oplus\ldots\oplus
E_{\lambda_d}.
\end{equation}
Since $\overset{.}w_i\cdot v_{\lambda_i}\in
E_{\lambda_i}$ is a weight vector  of
weight $w_i\cdot \lambda_i$, it follows
from (ii) and Lemma~\ref{closedness} that
the orbit $G\cdot v$ is closed in
$X_{\lambda_1,\ldots, \lambda_d}$. In turn,
this implies, by Matsushima's criterion,
see, e.g.,\,\cite[
Theorem\,4.17]{popov-vinberg2}, that $G_v$ is a
reductive group. We claim that $G_v$ is finite, i.e.,
that ${\rm Lie}(G_v)=0$.

To prove this, observe that since
$G_{\overset{.}w_i\cdot v_{\lambda_i}}=
\overset{.}w_i
G_{v_{\lambda_i}}{\overset{.}w}_i^{-1}$,
decomposition \eqref{point} implies that
\begin{equation}\label{decomp}
G_v=\bigcap_{i=1}^{d} \overset{.}w_i
G_{v_{\lambda_i}}{\overset{.}w}_i^{-1}.
\end{equation}
Taking into account that $\overset{.}w_i\in
N_G(T)$ and $G_{v_{\lambda_i}}$ is
normalized by $T$, we deduce from
\eqref{decomp} that $G_v$ is normalized by
$T$ as well. Hence, cf.,
e.g.,\,\cite[20.7]{TY},
\begin{equation}\label{Lie}
{\rm Lie}(G_v)=\mathfrak h\oplus \Bigl(
\underset{\alpha\in S}\oplus\mathfrak
g_{\alpha}\Bigr),
\end{equation}
where  $\mathfrak g_{\alpha}$ is the Lie
algebra of one-dimensional unipotent root
subgroup of $G$ corresponding to the root
$\alpha\in R$, $S$ is a subset of $R$, and
$\mathfrak h$ is a maximal torus of ${\rm
Lie}(G_v)$ contained in ${\rm Lie}(T)$.
Since  ${\rm Lie}(G_v)$ is reductive, the
conditions ${\rm Lie}(G_v)=0$ and
$\mathfrak h=0$ are equivalent. To prove
that $\mathfrak h=0$, observe that
\begin{equation}\label{decomp2}
{\rm Lie}(G_{v_{\lambda_i}})={\rm
Lie}\,({\rm ker}\,\lambda_i)  \oplus
\Bigl(\underset{\alpha\in
S_i}\oplus\mathfrak g_{\alpha}\Bigr)
\end{equation}
for some $S_i \subset R$, see
\cite{popov-vinberg1}. From \eqref{decomp}
and \eqref{decomp2} we then deduce that
\begin{equation}\label{conj}
{\rm Lie} (G_v)=\bigcap_{i=1}^{d}
\Bigl({\rm Lie}\,\bigl({\rm ker} (w_i\cdot
\lambda_i)\bigr)\oplus
\Bigl(\underset{\alpha\in w_i\cdot
S_i}\oplus\mathfrak g_\alpha\Bigr)\Bigr).
\end{equation}
In turn, it follows from \eqref{Lie} and
\eqref{conj} that
\begin{equation}\label{h}
\mathfrak h\subseteq \bigcap_{i=1}^{d} {\rm
Lie}\,\bigl({\rm
ker}(w_i\cdot\lambda_i)\bigr).
\end{equation}
From property (i) we deduce that the
right-hand side of \eqref{h} is equal to
$0$. Hence $\mathfrak h=0$, as claimed.
Thus we proved that $G_v$ is finite.

It follows from $\dim (G_v)=0$ that $\dim (G\cdot
v)=\dim (G)$. Hence maximum of dimensions of
$G$-orbits in $X_{\lambda_1,\ldots,\lambda_d}$ is
equal to $\dim (G)$. But the set of point whose
$G$-orbit has maximal dimension is open in
$X_{\lambda_1,\ldots,\lambda_d}$, cf.,
e.g.,\,\cite[1.4]{popov-vinberg2}. Hence
$G$-stabilizer of a point in general position in
$X_{\lambda_1,\ldots,\lambda_d}$ is finite. Finally,
since $G\cdot v$ is a closed orbit of maximal
dimension, \cite[
Theorem\,4]{popov0} implies that the
action of $G$ on $X_{\lambda_1,\ldots,\lambda_d}$ is
stable. \quad $\square$

\end{enavant}


\begin{enavant} {\bf Proof of Theorem
\ref{fin}}

Let $(\lambda_1,\ldots,\lambda_d)\in {\rm
P}_{\gg}^d$ be a primitive $d$-tuple.
Assume the contrary, i.e.,
\begin{equation}\label{>=}
{\rm sep}(G)+2\leqslant d. \end{equation}
From \eqref{>=} and Theorems \ref{teo1},
\ref{fntn} we deduce that the multiple flag
variety $G/P_{\lambda_1}\times\ldots\times
G/P_{\lambda_d}$ contains an open
$G$-orbit. Theorem \ref{bounds} then
implies that
\begin{equation}\label{=<}
d\leqslant {\rm rk}(G)+2.\end{equation}
From \eqref{>=}, \eqref{=<} we obtain the
inequality ${\rm sep}(G)\leqslant {\rm
rk}(G)$ that contradicts \eqref{sG}. \quad
$\square$
\end{enavant}


\begin{enavant}{\bf  Proof of
Theorem~{\bf\ref{concr}}}

By Theorem~\ref{simple} the claim follows
from the fact that  in either of the cases
listed in Table 2 the multiple flag variety
$G/P_{\lambda_1}\times
G/P_{\lambda_2}\times G/P_{\lambda_3}$
contains an open $G$-orbit. The latter is
proved as follows.

If $G$ is of type ${\sf B}_l$, ${\sf D}_l$,
${\sf E}_6$, or ${\sf E}_7$, then
$s_3=\{1,\ldots, {\rm rk}(G)\}$, hence
$P_{\lambda_3}=B$. Therefore
$G/P_{\lambda_1}\times
G/P_{\lambda_2}\times G/P_{\lambda_3}$
contains an open $G$-orbit if and only if
$G/P_{\lambda_1}\times G/P_{\lambda_2}$
contains an open $B$-orbit,
cf.,\,e.g.,\,\cite[Lem.\,4]{popov2}. All
the pairs of fundamental weights
$(\lambda_1,\lambda_2)$ for which the
latter holds are classified in
\cite[1.2]{littelmann}. According to this
classification, for these types of $G$, the
supports of $\lambda_1$ and $\lambda_2$ are
precisely (up to automorphism of the Dynkin
diagram) $s_1$ and $s_2$ specified in
Table~2.

For $G$ of types ${\sf A}_l$ and ${\sf
C}_l$, in \cite{MWZ1} and \cite{MWZ2} it is
given a classification of all the products
$G/P_{\lambda_1}\times
G/P_{\lambda_2}\times G/P_{\lambda_3}$ that
contain only finitely many $G$-orbits. One
of these orbits is then open in
$G/P_{\lambda_1}\times
G/P_{\lambda_2}\times G/P_{\lambda_3}$. The
triples $(\lambda_1, \lambda_2, \lambda_3)$
arising in these classifications are
precisely (up to automorphism of the Dynkin
diagram) the ones whose supports satisfy
the conditions of cases listed in Table 2
for these types of $G$. (Actually, in
\cite{MWZ1} and \cite{MWZ2}, flag varieties
are described in terms of ``compositions'',
i.e., essentially, dimension vectors of
corresponding flags. The information in
Table 2 is obtained by reformulating
results of \cite{MWZ1} and \cite{MWZ2} in
terms of supports of the corresponding
dominant weights; obtaining this
reformulation is not difficult: for
instance, for $G$ of type ${\sf A}_l$, one
deduces it from the fact that cardinality
of the set of nonzero parts of a
composition is equal to cardinality of the
support of corresponding dominant weight
plus $1$.) \quad $\square$
\end{enavant}


\begin{enavant} {\bf  Proof of
Theorem~\ref{geom}} \nopagebreak

Since $(0,\ldots,0)$ is a fixed point for
the action of $G$ on
$X_{\lambda_1,\ldots,\lambda_d}$, the
equivalence (iii)$\Leftrightarrow$(iv)
follows from the property that for every
reductive group action on affine variety,
disjoint invariant closed subsets are
separated by the algebra of invariants,
see,\,e.g.,\,\cite[
Theorem\,4.7]{popov-vinberg2}.

The equivalence (i)$\Leftrightarrow$(iv)
follows from \eqref{mgraduXC}, \eqref{XGG},
\eqref{**}, and Definition \ref{def2}.

The implication (i)$\Rightarrow$(ii)
 follows from Definition
\ref{def2}.

Arguing on the contrary, assume that (ii)
holds, but (i) does not. The latter means
that $c_{s_1\lambda_1,\ldots,
s_d\lambda_d}^0\geqslant 1$ for some
$(s_1,\ldots, s_d)\in \mathbb Z_{\geqslant
0}^d$. Lemma \ref{semigr} then implies that
$c_{ms_1,\ldots, ms_d}^0\geqslant 1$ for
every $m\in \mathbb Z_{>0}$. Taking
$m=m_1\ldots m_d$, we obtain
\begin{equation}\label{ccc}
c_{n_1 m_1\lambda_1,\ldots, n_d
m_d\lambda_d}^0\geqslant 1\quad\mbox{where
$n_i=m_1\ldots\widehat{m_i}\ldots m_ds_i$.}
\end{equation}
Definition \ref{def2} now shows that
property (ii) contradicts \eqref{ccc}.
\quad $\square$.
\end{enavant}

\begin{enavant} {\bf  Proof of
Theorem~\ref{bound}}

We utilize the following lemma.

\begin{lemma}\label{suter}
Let $R$ be an irreducible reduced root
system in an $n$-dimensional rational
vector space $L$. Then there are the Weyl
chambers $C_1,\ldots, C_{n+1}\subset L$ of
$R$ such that
\begin{equation}\label{in}
\mbox{$0\in {\rm conv}(\{x_1,\ldots,
x_{n+1}\})$
for every choice of points $x_1\in
\overline{C_1},\ldots, x_{n+1}\in
\overline{C_{n+1}}$.}
\end{equation}
\end{lemma}

\noindent {\it Proof}\;
 (R.\,Suter). Let $R^\vee\subset L^*$
 be the dual
 root system. Take a basis $l_1,\ldots,
 l_{n}$ of $R^\vee$ and let
 $-l_{n+1}$ be
 the corresponding maximal
 root of $R^\vee$.
 Utilizing notation \eqref{+-}, put
\begin{equation}\label{Z}
Z_i:=\bigcap_{j\in [1, n+1],\, j\neq i}
l_j^+.
\end{equation}
We claim that $0\in {\rm
conv}(\{x_1,\ldots, x_{n+1}\})$ for every
choice of points $x_i\in Z_i$, $i=1,\ldots,
n+1$. Indeed, if $0\notin {\rm
conv}(\{x_1,\ldots, x_{n+1}\})$ for some
$x_i\in Z_i$, $i=1,\ldots, n+1$, then there
is a nonzero linear function $l\in L^*$
such that
\begin{equation}\label{l<}
l(x_i)<0\quad \mbox{for every $i=1,\ldots,
n+1$}. \end{equation} Since
$L^*=\bigcup_{j=1}^{n+1}
 {\rm
 cone}(\{l_1,\ldots,\widehat{l_j},\ldots,
 l_{n+1}\})$,
there is $i_0$ such that
\begin{equation}\label{lco}
l\in {\rm
 cone}(\{l_1,\ldots,\widehat{l_{i_0}},\ldots,
 l_{n+1}\}).
\end{equation}
From \eqref{Z} and \eqref{lco} we deduce
that  $l(x_{i_0})\geqslant 0$, contrary to
\eqref{l<}. A contradiction.

Now, since every $Z_i$ is a union of the
closures of Weyl chambers, we can choose a
Weyl chamber $C_i$ lying in $Z_i$. Then
required property \eqref{in} holds for
$C_1, \ldots, C_{n+1}$.
 \quad $\square$

\vskip 2mm

Passing to the proof of Theorem
~\ref{bound} and arguing on the contrary,
assume that (i) fails, i.e., for some $i$,
\begin{equation}\label{con}
m_i\lambda_i^*=m_1\lambda_1+\ldots+\widehat{m_i\lambda_i}+\ldots+
m_d\lambda_d,
\end{equation}
where $(m_1,\ldots, m_d)\in \mathbb
Q_{\geqslant 0}^d$ and $m_i>0$. Multiplying
both sides of \eqref{con} by an appropriate
natural number, we may (and shall) assume
that $(m_1,\ldots, m_d)\in \mathbb
Z_{\geqslant 0}^d$. The Cartan
com\-po\-nent of
$E_{m_1\lambda_1}\otimes\ldots
\otimes\widehat{E_{m_i\lambda_i}}\otimes\ldots\otimes
E_{m_d\lambda_d}$ is
$E_{m_1\lambda_1+\ldots+
\widehat{m_i\lambda_i}+\ldots+
m_d\lambda_d}$. Hence
\begin{equation}\label{Cart}
E_{m_1\lambda_1}\otimes\ldots\otimes
\widehat{E_{m_i\lambda_i}}\otimes\ldots\otimes
E_{m_d\lambda_d}\simeq
E_{m_1\lambda_1+\ldots+\widehat{m_i\lambda_i}+\ldots+
m_d\lambda_d}\oplus\ldots,
\end{equation}
where the right-hand side of \eqref{Cart}
is a direct sum of simple $G$-modules. It
follows from \eqref{con}, \eqref{Cart}, and
\eqref{multipl}  that
\begin{equation}\label{ccontr}
c_{m_1\lambda_1,\ldots,
m_d\lambda_d}^0\geqslant 1.
\end{equation}
Since \eqref{ccontr} contradicts the
assumption that
$(\lambda_1,\ldots,\lambda_d)$ is
invariant-free, this proves~(i).

Again arguing on the contrary, assume that
(ii) fails, i.e.,
\begin{equation}\label{>rk}
d\geqslant r+1,\quad \mbox{where $r={\rm
rk}(G)$}.
\end{equation}
By Lemma \ref{suter} there are the Weyl
chambers $C_1,\ldots, C_{r+1}\subset
\mathcal X(T)_{\mathbb Q}$ of the root
system of $G$ with respect to $T$ such that
property \eqref{in} (with $n=r$) holds.
Inequality \eqref{>rk} implies that there
are (unique) elements $w_1,\ldots,
w_{r+1}\in W$ such that $w_i\cdot
\lambda_{i}\in \overline{C_i}$ for every
$i$. By \eqref{in} we have
\begin{equation*}\label{0in}
0\in {\rm conv}(\{w_1\cdot
\lambda_1,\ldots, w_{r+1}\cdot
\lambda_{r+1}\}).
\end{equation*}
Hence $0$ is an interior point of some face
of the potytope ${\rm
conv}(\{w_1\cdot\lambda_1,\ldots,
w_{r+1}\cdot\lambda_{r+1}\})$; whence
\begin{equation}\label{int0}
0\in {\rm int}\bigl({\rm
conv}(\{w_{i_1}\cdot\lambda_{i_1},\ldots,
w_{i_m}\cdot\lambda_{i_m}\})\bigr)
\end{equation}
for some $i_1,\ldots, i_m$. Since
$\overset{.}w_i\cdot v_{\lambda_i}\in
E_{\lambda_i}$ is a weight vector  of
weight $w_i\cdot\lambda_i$, it follows from
\eqref{int0} and Lemma~\ref{closedness}
that the $G$-orbit of point
$\overset{.}w_{i_1}\cdot v_{\lambda_{i_1}}+
\ldots + \overset{.}w_{i_m}\cdot
v_{\lambda_{i_m}}$ is closed in
$X_{\lambda_{i_1},\ldots, \lambda_{i_m}}$.
But $X_{\lambda_{i_1},\ldots,
\lambda_{i_m}}$ clearly admits a  closed
$G$-invariant embedding in
$X_{\lambda_1,\ldots, \lambda_d}$, so this
gives a closed $G$-orbit in
$X_{\lambda_1,\ldots, \lambda_d}$ as well.
Since this orbit is different from
$(0,\ldots, 0)$, Theorem \ref{geom} yields
a contradiction with the assumption that
$(\lambda_1,\ldots, \lambda_d)$ is
invariant-free. This proves (ii).\quad
$\square$

\end{enavant}


\begin{enavant}{\bf Stability of
\boldmath $G$-action on
$X_{\lambda_1,\ldots, \lambda_d}$}
\label{stabili}

In this section we prove that several other
conditions are sufficient for stability of
the action of $G$ on $X_{\lambda_1,\ldots,
\lambda_d}$.

Consider ${\rm P}_{++}^d$ as a submonoid of
the group $\mathcal X(T)^d$ that, in turn,
is considered as a lattice in the rational
vector space $\mathcal X(T)^d_{\mathbb
Q}:=\mathcal X(T)^d\otimes\mathbb Q$.
Notice that if $A$ and $B$ are submonoids
of ${\rm P}_{++}^d$, then the condition
\begin{equation*}
{\rm int}\bigl({\rm cone}(A)\bigr)\cap {\rm
int}\bigl({\rm cone}(B)\bigr)\neq
\varnothing
\end{equation*}
is equivalent to the property that
$A-B:=\{a-b\mid a\in A, b\in B\}$ is a
group.

 For
$(\lambda_1,\ldots,\lambda_d)\in{\rm
P}_{++}^d$, consider the submonoid
$\langle\lambda_1,\ldots,\lambda_d\rangle$
of ${\rm P}_{++}^d$ generated by
$(\lambda_1,0,\ldots, 0),\ldots, (0,\ldots,
0, \lambda_d)$,
$$
\langle\lambda_1,\ldots,\lambda_d\rangle:=
\{(n_1\lambda_1, \ldots, n_d\lambda_d)\mid
(n_1,\ldots, n_d)\in \mathbb Z_{\geqslant 0}^d\}.
$$
Put
\begin{equation}\label{1111}
\Gamma (G,d):=\{(\mu_1,\ldots,\mu_d)\in{\rm
P}_{++}^d\mid (E_{\mu_1}\otimes
\ldots\otimes E_{\mu_d})\hskip -.1mm^G \neq
0 \}.
\end{equation}
\begin{example}\label{ggm}
$\Gamma(G,1)=\{0\}$, and, by
\eqref{multipl}, we have $
\Gamma(G,2)=\{(\mu, \mu^*)\mid \mu\in {\rm
P}_{++}\}.$
\end{example}

\begin{example} Put ${\rm LR}(G,
3):=\{(\lambda_1,\lambda_2,\lambda_3)\mid
(\lambda_1,\lambda_2,\lambda_3^*)\in \Gamma(G, 3)\}$.
Then ${\rm LR}({\bf SL}_n, 3)$ is the
Littelwood--Richardson semigroup of order $n$,
\cite{zel}. It has been intensively studied during the
last decade and is now rather well understood. For
instance, a minimal system of linear inequalities
cutting out ${\rm cone}\bigl({\rm LR}({\bf SL}_n,
3)\bigr)$ in $\mathcal X(T)_\mathbb Q^3$ is found, the
walls of ${\rm cone}\bigl({\rm LR}({\bf SL}_n, 3)\bigr)$
are described, and it is proved that ${\rm LR}({\bf
SL}_n, 3)$ is the intersection of ${\rm cone}\bigl({\rm
LR}({\bf SL}_n, 3)\bigr)$ with the corresponding lattice
in $\mathcal X(T)_\mathbb Q^3$ (saturation conjecture),
see survey \cite{fulton} and \cite{belkale1},
\cite{belkale2}. This immediately implies analogous
results about $\Gamma({\bf SL}_n, 3)$. In \cite{km} some
general structural results for $\Gamma(G, 3)$ are
obtained and $\Gamma({\bf Sp}_4,3)$ and $\Gamma({\bf
G}_2,3)$ are computed. $\Gamma({\bf Spin}_8, 3)$ is
studied in \cite{kkm}.
\end{example}

These examples show that $\Gamma (G, d)$
for $d\leqslant 3$ is a finitely generated
submonoid of ${\rm P}_{++}^d$. Actually
this is true for every $d$, see Corollary
of Theorem \ref{GA} below. It would be
interesting to understand the structure of
this monoid in general case. What are the
inequalities cutting out ${\rm
cone}\bigl(\Gamma\,(G, d)\bigr)$ in
$\mathcal X(T)_\mathbb Q^d$? What are the
generators of $\Gamma (G, d)$?

\begin{theorem} \label{GA} Consider
$G$ as the diagonal subgroup of $G^d$. Then
{\rm(}see {\rm \eqref{monoid}}{\rm)}
\begin{equation}\label{gamma}
\Gamma(G, d)=\mathcal S(G^d, G^d/G).
\end{equation}
\end{theorem}

\begin{proof} We can (and shall)
identify in the natural way ${\rm
P}_{++}^d$ with the monoid of dominant
weights of the semisimple group $G^d$ with
respect to maximal torus $T^d$ and Borel
subgroup $B^d$. Simple $G^d$-modules are
tensor products
$E_{\mu_1}\otimes\ldots\otimes E_{\mu_d}$,
where $E_{\mu_i}$ is considered as the
$G^d$-module via the $i$th projection
$G^d\rightarrow G$, cf.,
e.g.,\,\cite[Ch.\,4, \S3]{OV}. This,
Frobenius duality
(cf.,\,e.g.,\,\cite[
Theorem\,3.12] {popov-vinberg2}), and formulas
\eqref{**}, \eqref{1111}, \eqref{monoid} now imply
the claim. \quad $\square$
\renewcommand{\qed}{}
\end{proof}

\noindent{\bf Corollary.} {\it
$\Gamma(G,d)$ is a finitely generated
submonoid of ${\rm P}_{++}^d$}.

\begin{proof} Equality \eqref{gamma}
implies
that $\Gamma(G,d)$  is a
 submonoid of ${\rm
P}_{++}^d$. Since $G$ is a reductive group,
Matsushima's criterion implies that $G^d/G$
is an affine variety; whence $\Gamma(G,d)$
is finitely generated (see the arguments
right after formula \eqref{monoid}). \quad
$\square$
\renewcommand{\qed}{}
\end{proof}

\begin{theorem}\label{>>>>} Let $(\lambda_1,\ldots,
\lambda_d)\in {\rm P}_{\gg}^d$. If either
of the following conditions holds, then the
action of $G$ on
$X_{\lambda_1,\ldots,\lambda_d}$ is stable:
\begin{enumerate}
\item[\rm(i)] ${\rm int}\bigl({\rm
cone}\bigl(\Gamma(G,d)\bigr)\bigr)\cap {\rm
int}\bigl({\rm
cone}\bigl(\{(\lambda_1,0,\ldots,
0),\ldots,
(0,\ldots,0,\lambda_d)\}\bigr)\bigr)\neq
\varnothing$; \item[\rm(ii)] there is $i$
such that
\begin{gather}\label{x}
\begin{gathered} \dim\bigl({\rm
cone}(\{\lambda_1,\ldots,
\widehat{\lambda_i},\ldots,
\lambda_d\})\bigr)={\rm rk}(G),\\
\lambda_i^*\in {\rm int}\bigl({\rm
cone}(\{\lambda_1,\ldots,
\widehat{\lambda_i},\ldots,
\lambda_d\})\bigr);
\end{gathered}
\end{gather}
\item[\rm(iii)] $\{1,\ldots, d\}$ is a
disjoint union of subsets $\{i_1,\ldots,
i_s\}$ and $\{j_1,\ldots, j_t\}$ such that
\begin{gather}\label{xx}
\begin{gathered}
\dim\bigl({\rm
cone}(\{\lambda_{i_1},\ldots,
\lambda_{i_s}\})\bigr)= \dim\bigl({\rm
cone}(\{\lambda_{j_1},\ldots,
\lambda_{j_t}\})\bigr)={\rm rk}(G),\\ {\rm
int}\bigl({\rm
cone}(\{\lambda_{i_1},\ldots,
\lambda_{i_s}\})\bigr)\cap {\rm
int}\bigl({\rm
cone}(\{\lambda_{j_1}^*,\ldots,
\lambda_{j_t}^*\})\bigr)\neq\varnothing.
\end{gathered}
\end{gather}
\end{enumerate}
\end{theorem}
\begin{proof}
(1) Discussion in Section \ref{po} (see
formula \eqref{isomo}) implies that
$k[X_{\lambda_1,\ldots,
\lambda_d}]_{(n_1,\ldots, n_d)}$ is a
simple $G^d$-module with highest weight
$(n_1\lambda_1^*,\ldots,n_d\lambda_d^*)$.
This and \eqref{mgradu} imply
 that
\begin{equation}\label{sd} \mathcal S(G^d,
X_{\lambda_1,\ldots,
\lambda_d})=\langle\lambda_1^*,\ldots,
\lambda_d^*\rangle.\end{equation} By
\cite[
Theorem\,10]{V}  the action of $G$ on
$X_{\lambda_1,\ldots,\lambda_d}$ is stable if
$\mathcal S(G^d,X_{\lambda_1,\ldots,\lambda_d})-
\mathcal S(G^d, G^d/G)$ is a group. But $\Gamma(G,
d)^*=\Gamma(G, d)$ by \eqref{**} and \eqref{1111}.
Hence \eqref{sd} and Theorem \ref{GA} imply that the
action of $G$ on $X_{\lambda_1,\ldots,\lambda_d}$ is
stable if (i) holds.

(2)  The variety $X_{\lambda_1,\ldots,
\lambda_d}$ is $G$-isomorphic to $Y\times
Z$, where $Y:=X_{\lambda_1,\ldots,
\widehat{\lambda_i},\ldots, \lambda_d}$ and
$Z:=X_{\lambda_i}$. Discussion in Section
\ref{po} shows that $\mathcal S(G,
Y)^*\ni\lambda_1,\ldots,
\widehat{\lambda_i},\ldots, \lambda_d$ and
$\mathcal S(G, Z)=\mathbb Z_{\geqslant
0}\lambda_i^*$. Hence
\begin{equation}\label{cap}{\rm cone}
\bigl(\mathcal S(G, Y)^*\bigr)\supseteq
{\rm cone}(\{\lambda_1,\ldots,
\widehat{\lambda_i},\ldots,
\lambda_d\})\quad \mbox{and}\quad {\rm
cone}\bigl(\mathcal S(G, Z)\bigr)=\mathbb
Q_{\geqslant 0}\lambda_i^*.\end{equation}
If (ii) holds, we deduce from \eqref{x} and
\eqref{cap} that
\begin{equation}\label{capcap}
{\rm int}\bigl({\rm cone}\bigl(\mathcal
S(G, Y)^*\bigr)\bigr)\cap {\rm
int}\bigl({\rm cone}\bigl(\mathcal S(G,
Z)\bigr)\bigr)\neq
\varnothing.\end{equation} By
\cite[
Theorem\,9]{V} inequality \eqref{capcap} implies that
the action of $G$ on $Y\times Z$ is stable.

(3) Assume now that (iii) holds. The
variety $X_{\lambda_1,\ldots, \lambda_d}$
is isomorphic to $Y\times Z$, where $Y:=
X_{\lambda_{i_1},\ldots, \lambda_{i_s}}$
and $Z:=X_{\lambda_{j_1},\ldots,
\lambda_{j_t}}$. Hence
\begin{gather}\label{capcapcap}
\begin{gathered}{\rm cone}\bigl(\mathcal S(G,
Y)^*\bigr)\supseteq {\rm
cone}(\{\lambda_{i_1},\ldots,
\lambda_{i_s}\}),\\
{\rm cone}\bigl(\mathcal S(G,
Z)\bigr)\supseteq {\rm
cone}(\{\lambda_{j_1}^*,\ldots,
\lambda_{j_t}^*\}).
\end{gathered}
\end{gather}
It follows from \eqref{xx} and
\eqref{capcapcap} that, as above,
\eqref{capcap} holds and hence the action
of $G$ on $Y\times Z$ is stable. \quad
$\square$
\renewcommand{\qed}{}
\end{proof}

\end{enavant}


\begin{enavant}{\bf Case
of \boldmath${\bf SL}_{n}$}\label{12}

Let $G={\bf SL}_{n}$. In this case,
combining the above results with that of
the re\-presentation theory of quivers (we
refer to \cite{kac1}, \cite{kac2},
\cite{derksen-weyman}, \cite{sch1} for the
notions of this theory) leads to a
characterization of primitive $d$-tuples of
fundamental weights in terms of canonical
decomposition of dimension vectors of
representations of some graphs and to an
algorithmic way of solving, for every such
$d$-tuple $(\varpi_{i_1},\ldots,
\varpi_{i_d})$, whether it is primitive  or
not.

Namely, in this case, $G/P_{\varpi_i}$ is
the Grassmannian variety of $i$-dimensional
linear subspaces on $k^{n}$, and the
existence of an open $G$-orbit in $
 G/P_{\varpi_{i_1}}
 \times
 \ldots\times G/P_{\varpi_{i_d}}
 $
admits the following reformulation in terms
of the representation theory of quivers.
Let ${\mathcal V}_d$ be the quiver with
$d+1$ vertices, $d$ outside, one inside,
and the arrows from each vertex outside to
a vertex inside (the vertices are
enumerated by $1,\ldots, d+1$ so that the
inside vertex is enumerated by $1$):
\begin{center}
\leavevmode \epsfxsize =3cm
\epsffile{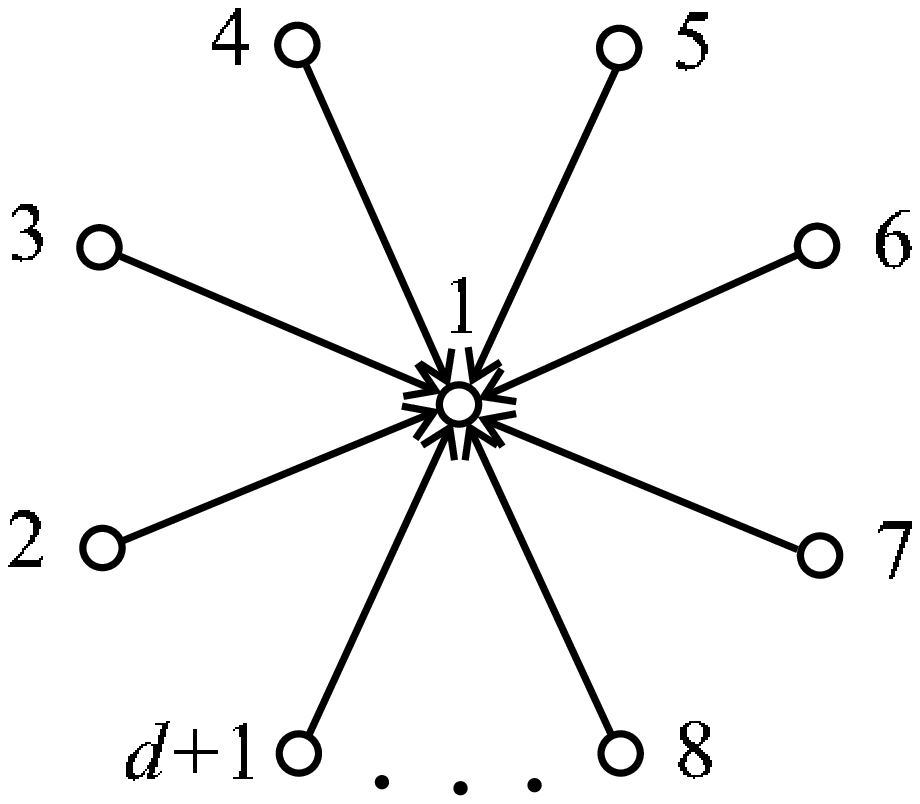}
\end{center}
\noindent
 Given a vector
\begin{equation*}
\alpha:=(a_1,\ldots, a_{d+1})\in \mathbb
Z_{\geqslant 0}^{d+1},
\end{equation*}
put ${\bf GL}_{\alpha}: ={\bf GL}_{a_1}
\times\ldots\times {\bf GL}_{a_{d+1}} $ (we
set ${\bf GL}_0:=\{e\}$). Let
\begin{equation*}{\rm Rep}(\mathcal
V_d, \alpha):= {\rm Mat}_{a_1\times a_2}
\oplus\ldots\oplus{\rm Mat}_{a_1\times
a_{d+1}}
\end{equation*} be
the space of $\alpha$-dimensional
representations of ${\mathcal V}_d$ endowed
with the natural ${\bf
GL}_{\alpha}$-action. For $\mathcal V_d$,
the Euler inner product $\langle\ {\,|\,}\
\rangle$ on $\mathbb Z^{d+1}$ is given by
\begin{equation}\label{euler}
\langle(x_1,\ldots
x_{d+1})\,\vert\,(y_1,\ldots,
y_{d+1})\rangle=(x_1y_1+\ldots
+x_{d+1}y_{d+1})-y_1(x_2+\ldots +x_{d+1}).
\end{equation}

It then follows from the basic definitions
that the following properties are
equi\-va\-lent:
 \begin{enumerate}
 \item[(a)] the multiple flag variety
 $
 G/P_{\varpi_{i_1}}
 \times
 \ldots\times G/P_{\varpi_{i_d}}
 $ contains an open $G$-orbit;
 \item[(b)] the space
 ${\rm
Rep}(\mathcal V_d, \gamma)$, where
$$\gamma:=(n,
 i_1,\ldots, i_d),$$  contains
an open ${\bf GL}_\gamma$-orbit.
\end{enumerate}
\begin{theorem}\label{sl-fund} Let $G={\bf
SL}_{n}$. The following properties are
equivalent:
\begin{enumerate}
\item[(i)]
$(\varpi_{i_1},\ldots,\varpi_{i_d})$ is
primitive; \item[(ii)] all the roots
$\beta_i$ appearing in the canonical
decomposition of $\gamma$,
\begin{equation}\label{canon}
\gamma=\beta_1+\ldots+\beta_s,\end{equation}
\noindent are real, i.e.,
$\langle\beta_i\,|\,\beta_i\rangle=1$.
\end{enumerate}
\end{theorem}
\begin{proof}
By Theorems~\ref{teo1} and \ref{teo2}
properties (i) and (a) are equivalent. On
the other hand, properties (ii) and (b) are
equivalent by \cite[Cor.\,1 of
Prop.\,4]{kac2}.
 \quad
$\square$
\renewcommand{\qed}{}
\end{proof}

Note that there are combinatorial
algorithms for finding decomposition
\eqref{canon} (see \cite{sch1},
\cite{derksen-weyman}; the algorithm in
\cite{derksen-weyman} is fast). Hence they,
Theorem~\ref{sl-fund}, and formula
 \eqref{euler}
 yield algorithms
 verifying, for
every
concrete $d$-tuple
$(\varpi_{i_1},\ldots,\varpi_{i_d})$, whether it is
primitive  or~not.
\end{enavant}


\renewcommand{\baselinestretch}{1.1}

 \providecommand{\bysame}{\leavevmode\hbox
to3em{\hrulefill}\thinspace}


\begin{thebibliography}{XXXXx}

\bibitem[Be${}_1$]{belkale1}{\sc P.\,Belkale},
\emph{Geometric proof of a conjecture of Fulton}, {\tt
arXiv:math.AG/0511664.}

\bibitem[Be${}_2$]{belkale2}{\sc P.\,Belkale},
\emph{Geometric proofs of Horn and saturation
conjectures}, {\tt arXiv:math.AG/0208107.}

\bibitem[Bo${}_1$]{bourbaki1}
  {\sc  N.\,Bourbaki},
  \emph{Alg\`ebre commutative}, Chap. 5,
  Hermann, Paris, 1964.


  \bibitem[Bo${}_2$]{bourbaki}
  {\sc  N.\,Bourbaki},
  \emph{Groupes et alg\`ebres de Lie},
  Chap. IV, V,  VI, Hermann, Paris, 1968.

  \bibitem[DW]{derksen-weyman}
  {\sc H.\,Derksen, J.\,Weyman},
  \emph{On the canonical decomposition of quiver
  representations},
  Compositio Math. {\bf 133} (2002), 245--265.

\bibitem[F]{fulton} {\sc W.\,Fulton},
\emph{Eigenvalues, invariant factors,
highest weights, and Schubert calculus},
Bull. Amer. Math. Soc. {\bf 37} (2000), No.
3, 209--249.


  \bibitem[K${}_1$]{kac1}
  {\sc V.\,Kac},
  \emph{Infinite root systems,
  representations of graphs and invariant
  theory}, Invent. Math. {\bf 56}
  (1980), 57--92.

  \bibitem[K${}_2$]{kac2}
  {\sc V.\,Kac},
  \emph{Infinite root systems,
  representations of graphs and invariant
  theory} II, J. Algebra {\bf 78}
  (1982), 141--162.

  \bibitem[KM]{km}
  {\sc M.\,Kapovich, J.\,Millson},
  \emph{Structure of the
  tensor product semigroup}, {\tt
  arXiv:math.RT/050}
  {\tt 8186.}

  \bibitem[KKM]{kkm}
  {\sc M.\,Kapovich, S.\,Kumar, J.\,Millson},
  \emph{Saturation and irredundancy for ${\rm Spin}(8)$},
  {\tt
  arXiv:math.}
  {\tt RT/0607454}.


  \bibitem[Li]{littelmann}
{\sc P.\,Littelmann}, \emph{On spherical
double cones}, J. Algebra {\bf 166} (1994),
no. 1, 142--157.

\bibitem[Lu]{luna}
{\sc  D.~Luna}, \emph{Slices \'etales},
Bull. Soc. Math. France, M\'emoire {\bf 33}
(1973), 81--105.

\bibitem[MWZ${}_1$]{MWZ1}
{\sc P.\,Magyar, J.\,Weyman,
A.~Zelevinsky}, \emph{Multiple flag
varieties of finite type}, Adv. Math. {\bf
141} (1999), no. 1, 97--118.

\bibitem[MWZ${}_2$]{MWZ2}
{\sc P.\,Magyar, J.\,Weyman,
A.~Zelevinsky}, \emph{Symplectic multiple
flag varieties of finite type}, J. Algebra
{\bf 230} (2000), no. 1, 245--265.

\bibitem[M]{M}
{\sc H.\,Matsumura}, \emph{Commutative
Algebra}, W.\,A.\,Benjamin Co., New York,
1970.

\bibitem[OV]{OV}
{\sc A.\,L.\,Onishchik, E.\,B.\,Vinberg},
\emph{Lie Groups and Algebraic Groups},
Springer-Verlag, Berlin, Heidelberg, 1990.

\bibitem[P${}_1$]{popov0}
{\sc V.\,L.\,Popov}, \emph{Stability
criteria for the action of a semisimple
group on a factorial manifold}, Math. USSR
Izv. {\bf 4} (1970), 527--535.


\bibitem[P${}_2$]{popov1}
{\sc V.\,L.\,Popov}, \emph{Picard groups of
homogeneous spaces of linear algebraic
groups and one-dimensional homogeneous
vector bundles}, Math. USSR, Izv. {\bf 8}
(1975), 301--327.

\bibitem[P${}_3$]{popov3}
{\sc V.\,L.\,Popov}, {\it On the closedness
of some orbits of algebraic groups}, Funct.
Anal. Appl. {\bf 31} (1997), No.\,4,
286-289.




  \bibitem[P${}_4$]{popov2}
  {\sc V.\,L.\,Popov}, \emph{Generically multiple
transitive algebraic group actions}, in: {\it Proc.
International Colloquium {\rm``}Algebraic Groups and
Homogeneous Spaces{\rm''}, $6$--$14$ January $2004$},
Tata Inst. Fund. Research, Bombay, India, Narosa
Publ. House, 2006, pp. 481--523.

\bibitem[PV${}_1$]{popov-vinberg1}
 {\sc  V.\,L.\,Popov,  E.\,B.\,Vinberg},
  \emph{On a class of quasihomogeneous
  affine
  varieties},
  Math. USSR, Izv. {\bf 6} (1973), 743--758.




  \bibitem[PV${}_2$]{popov-vinberg2} {\sc V.\,L.\,Popov,
E.\,B.\,Vinberg},
  \emph{Invariant Theory}, Encycl.
  Math. Sci., Vol.
  55, Springer-Verlag, Heidelberg, 1994,
  pp. 123--284.

\bibitem[R${}_1$]{rosenlicht0}
{\sc M.\,Rosenlicht}, \emph{Some
rationality questions on algebraic groups},
Ann. Mat. Pura Appl. (4) {\bf 43} (1957),
25--50.

\bibitem[R${}_2$] {rosenlicht} {\sc M.\,Rosenlicht},
\emph{A remark on quotient spaces}, An.
Acad. Brasil. Ci. {\bf 35} (1963),
487--489.

\bibitem[Sc]
{sch1}
  {\sc  A.\,Schofield}, \emph{General representations
  of quivers}, Proc. London Math. Soc.
  {\bf 65} (1992),
  46--64.

  \bibitem[Se]{serre} {\sc J.-P. Serre},
  \emph{Alg\`ebres de Lie semi-simples
  complexes}, Bejamin, New York, Amsterdam,
  1966.

\bibitem[St]{St}
{\sc J.\,R.\,Stembridge},
\emph{Multiplicity free products of Schur
functions}, Ann. Comb. {\bf 5} (2001),
no.\,2, 113--121.

\bibitem[TY]{TY}
{\sc P.\,Tauvel, R.\,W.\,T.\,Yu}, \emph{Lie
Algebras and Algebraic Groups}, Springer,
Berlin, Heidelberg, New York, 2005.

\bibitem[V]{V}
{\sc E.\,B.\,Vinberg}, \emph{On stability
of actions of reductive algebraic groups},
in: {\sc Fong, Yuen} (ed.) et al.,
\emph{Lie Algebras, Rings and Related
Topics}, Papers of the 2nd Tainan--Moscow
international algebra workshop'97, Tainan,
Taiwan, January 11--17, 1997, Hong Kong,
Springer, 2000, pp.\,188--202.



\bibitem[Z]{zel}
{\sc A.\,Zelevinsky},
\emph{Littelwood--Richardson semigroups},
in: {\it New Perspectives in Algebraic
Combinatorics}, Cambridge University Press
(MSRI Publication), 1999, pp.\,337--345,
{\tt arXiv:math. CO/9704228}.



\end{thebibliography}
\end{document}